\documentclass[10pt,reqno]{amsart}

\usepackage{amsmath,amssymb,amsfonts,mathtools,bm}
\usepackage{graphicx,version}
\graphicspath{{figures/}}
\usepackage{booktabs,tabularx,array}
\usepackage{enumitem}
\usepackage{xcolor}
\usepackage{microtype}
\usepackage{placeins}
\usepackage{url,hyperref}
\hypersetup{hidelinks}

\numberwithin{equation}{section}
\newtheorem{theorem}{Theorem}[section]
\newtheorem{lemma}[theorem]{Lemma}
\newtheorem{proposition}[theorem]{Proposition}
\newtheorem{corollary}[theorem]{Corollary}
\theoremstyle{definition}

\theoremstyle{remark}
\newtheorem{remark}[theorem]{Remark}

\newcommand{\R}{\mathbb{R}}
\newcommand{\cE}{\mathcal{E}}
\newcommand{\cK}{\mathcal{K}}

\newcommand{\dt}{\delta t}

\newcommand{\bu}{\bm{u}}
\newcommand{\bv}{\bm{v}}
\newcommand{\bfv}{\bm{f}}
\newcommand{\bg}{\bm{g}}
\newcommand{\bzero}{\bm{0}}
\newcommand{\curl}{\nabla\!\times}
\newcommand{\norm}[1]{\left\lVert #1\right\rVert}
\newcommand{\ip}[2]{\left(#1,#2\right)}

\begin{document}
\raggedbottom

\title[Long-time stable schemes for the NSE]
{High-order long-time stable methods for the forced incompressible Navier--Stokes equations}

\author[F. Huang, J. Shen and X. Wang]{Fukeng Huang$^\dagger$, Jie Shen$^\dagger$ and Xiaoming Wang$^\dagger$}
\thanks{$^\dagger$ School of Mathematical Sciences, Eastern Institute of Technology, Ningbo, Zhejiang 315200, P. R. China (fkhuang@eitech.edu.cn, jshen@eitech.edu.cn, wxm@eitech.edu.cn). This work was supported in part by National Natural Science Foundation of China (Grant Nos.
12501554, W2431008, 12531015 and 12371409).}


\begin{abstract}
We present in this paper two scalar auxiliary variable (SAV) correction
steps for constructing long-time stable schemes for the forced
incompressible Navier--Stokes equations with no-slip boundary conditions.
A main feature is that the same correction steps apply to any linear
predictor satisfying the order and spatial regularity requirements in
Assumption~P. Their construction and stability analysis are independent
of the particular form of the predictor and require no energy estimate
for it. The first correction uses a damped, shifted kinetic-energy
variable to obtain a uniform-in-time $L^2$ bound in two and three
dimensions, together with a weighted dissipation estimate. In two
dimensions, the second correction combines a higher-energy auxiliary
variable with a scalar relaxation to yield a uniform-in-time $H^1$
bound. Under the stated assumptions, these estimates hold on arbitrary
positive time grids, without restrictions on step sizes or adjacent
step ratios. We also derive eventual absorbing bounds; for the two-SAV
scheme, these apply to the admissible initial class associated with
fixed parameters. The generalized BDF schemes require one
constant-coefficient Helmholtz solve per velocity component and one
pressure Poisson solve per step; the second correction adds no field
solve. Numerical experiments exhibit second-, third-, and fourth-order
convergence for the tested BDF schemes and illustrate the long-time
stability and efficiency of adaptive time stepping.
\end{abstract}

\maketitle

\section{Introduction}
We consider in this paper the construction and long-time stability
analysis of efficient numerical schemes for the forced incompressible
Navier--Stokes equations
\begin{subequations}\label{eq:NS}
\begin{align}
 \bu_t-\nu\Delta\bu+(\bu\cdot\nabla)\bu+\nabla p&=\bfv,
    &&\text{in }\Omega\times(0,\infty),\label{eq:NS-mom}\\
 \nabla\cdot\bu&=0,
    &&\text{in }\Omega\times(0,\infty).\label{eq:NS-div}
\end{align}
\end{subequations}
Here $\bu$ and $p$ denote the velocity and pressure, respectively,
$\nu>0$ is the kinematic viscosity, and $\bfv$ is an external force.
We prescribe an initial velocity on a bounded connected domain
$\Omega\subset\R^d$, $d=2,3$, with a sufficiently smooth boundary.
For definiteness, we consider the no-slip boundary condition
$\bu=\bzero$ on $\partial\Omega$. The stability and absorption results
also extend to the mean-zero periodic setting described in
Section~\ref{sec:prelim}. Throughout the paper, we assume
\begin{equation}\label{eq:forcing}
 F_\infty:=\sup_{t\ge0}\norm{\bfv(\cdot,t)}_{L^2(\Omega)}<\infty.
\end{equation}

For long-time simulations, it is desirable that a numerical scheme
preserve the dissipative properties of the underlying equations. In
the presence of a persistent force, the kinetic energy need not
decrease, but the balance between forcing and viscous dissipation
still yields an $L^2$ absorbing bound in two and three dimensions and
an $H^1$ absorbing bound in two dimensions
\cite{ConstantinFoias,Temam,MajdaWang}. Unlike an all-time bound, which
may depend on the initial state, an absorbing bound gives a ball that
eventually contains the solution and whose radius is independent of
that state \cite{HumphriesStuart1994}. Our aim is to construct linear,
fully decoupled schemes that retain these properties and can be
combined with high-order time discretizations.

Long-time stability of Navier--Stokes discretizations has been studied
extensively. Early results include the convergence analysis of a
linearized backward Euler scheme for the unforced two-dimensional
problem \cite{Geveci}, and uniform bounds and discrete-attractor
convergence for a forced linear backward Euler scheme \cite{Ju}.
Other developments include implicit Euler and Crank--Nicolson methods
\cite{ToneEuler,ToneCN}, direct and fractional-step algorithms
\cite{SimoArmero}, and nonlinear Galerkin and Runge--Kutta
approximations \cite{JauberteauRosierTemam1990,Shen1990,LubichOstermann}.

For linearly extrapolated BDF2 schemes, Akbas, Kaya, and Rebholz
\cite{AkbasKayaRebholz} established all-time $L^2$ estimates, and Shiue
\cite{Shiue2025} obtained a refined unconditional long-time $L^2$
bound for the two-dimensional problem. Rebholz and Tone
\cite{RebholzTone} proved long-time $H^1$ stability in the
velocity--pressure formulation under a mild time-step restriction.
Unconditional long-time bounds were also obtained for first-order
and IMEX--BDF2 finite-element schemes in the velocity--vorticity
formulation \cite{Heister}. For the periodic vorticity equation,
semi-implicit multistep methods of order up to four allow fully
explicit convection and constant-coefficient solves, but their
long-time estimates require data-dependent time-step restrictions
\cite{GottliebToneWang,Wang2012,ChengWang}.

The scalar auxiliary variable (SAV) approach, introduced for gradient
flows \cite{ShenXuYangJCP,ShenXuYangSIREV}, provides an effective way
to treat nonlinear terms explicitly while retaining discrete energy
control. SAV methods for incompressible flows include the
auxiliary-energy-variable method of Lin, Yang, and Dong
\cite{LinYangDong}, which requires a nonlinear scalar solve, and the
linear, decoupled pressure-correction methods of Li, Shen, and Liu
\cite{LiShenLiu}. High-order SAV schemes have been developed for
periodic Navier--Stokes equations \cite{HuangShen2021}, nonperiodic
problems, and rotating flows \cite{wu2022new}. However, unconditional
stability on a fixed time interval does not imply uniform-in-time
boundedness for a forced problem. For example, the weak stability
estimate in \cite[Theorem~6]{HuangShen2023} and the higher-order error
estimates in \cite[Theorem~4.1]{HuangShenMathComp2026} involve
constants that depend on the final time.

Recent SAV methods address the forced problem by introducing
auxiliary damping. Han and Wang \cite{HanWang2024} constructed a
SAV--BDF2 scheme with uniform-in-time velocity $L^2$ bounds in two
and three dimensions and higher vorticity estimates in a
two-dimensional free-slip setting. Lin and Shiue
\cite{LinShiue2025} refined the $L^2$ analysis to obtain a
time-step-independent energy bound and to allow nonsmooth initial
data through regularization. For the two-dimensional periodic
vorticity equation, Wang, Wang, and Zhang
\cite{WangETD2026,WangLinearETD2026} developed second-order
variable-step exponential time-differencing methods with a
mean-reverting SAV. These methods require either a scalar cubic
solve or a linear auxiliary recovery by numerical inverse Laplace
transformation. We seek to combine auxiliary damping with high-order
predictors and decoupled constant-coefficient solves, while imposing
no time-step restriction in the long-time velocity estimates.

There are two main difficulties. The first is to obtain an $H^1$
bound for the no-slip problem. The provisional velocity is not
exactly divergence free, and the pressure and convection terms must
be controlled without the cancellations available in the periodic
vorticity formulation. The second is to establish a uniform estimate
on arbitrary positive time grids. On a nonuniform grid, shifting an
index in a time average changes its weights, while large individual
steps prevent the direct use of fixed-length discrete windows.

A main contribution of this paper is the construction of two scalar
correction steps whose formulation and long-time stability analysis
are independent of the particular linear predictor. The same
corrections can be combined with any linear predictor satisfying
Assumption~P. Thus, the generalized BDF schemes presented below are
particular realizations of a construction that also applies to other
multistep and Runge--Kutta predictors satisfying this assumption.

We first introduce a damped, shifted kinetic-energy variable. Its
scalar recursion gives an $L^2$ bound and a weighted dissipation
estimate without using an energy estimate for the predictor. In
two dimensions, we use the Stokes-pressure estimate of Liu, Liu, and
Pego \cite{LiuLiuPego} and interpolation inequalities to construct a
second scalar equation with a nonnegative residual. A scalar
relaxation transfers the first dissipation estimate to the second
auxiliary variable. We then apply a uniform Gronwall argument with
right-endpoint weights and treat large steps directly. This yields
a uniform $H^1$ bound without an upper time-step bound or an adjacent
step-ratio condition.

The main contributions of this paper are as follows.
\begin{enumerate}[label=(\roman*),leftmargin=*]
 \item We construct two SAV correction steps that can be combined
 with any linear predictor satisfying Assumption~P. Under the stated
 hypotheses, their stability estimates hold on arbitrary positive
 time grids and require no energy estimate for the predictor. The
 scalar analysis applies to any fixed correction degree.
 \item For the correction step with   one-SAV, we prove
 uniform-in-time velocity $L^2$ bounds in two and three dimensions.
 The estimate separates the decaying initial-data contribution from
 the forcing-dependent bound and yields an eventual $L^2$ absorbing
 radius (Theorem~\ref{thm:L2}).
 \item For the correction step with   two-SAV,  we prove a
 uniform-in-time velocity $H^1$ bound in two dimensions. We also
 establish eventual absorbing bounds for fixed algorithmic
 parameters and their admissible initial class
 (Theorem~\ref{thm:H1} and Corollary~\ref{cor:H1-absorption}).
\end{enumerate}
The schemes are easy to implement. The generalized BDF schemes require
one constant-coefficient Helmholtz solve for each velocity component
and one pressure Poisson solve per step. The two-SAV extension adds
only scalar operations and a velocity rescaling. Numerical tests
exhibit high-order accuracy and illustrate the effect of the
corrections on long-time stability and the efficiency of adaptive
time stepping. The convergence order of each coupled scheme requires
a separate analysis and depends on the predictor, the starting values,
and the consistency of the SAV ratios.

The rest of this paper is organized as follows.
Section~\ref{sec:prelim} introduces the notation, energy balances,
and linear predictor. We present and analyze the one-SAV and
two-SAV schemes in Sections~\ref{sec:one-sav}
and~\ref{sec:two-sav}, respectively. Section~\ref{sec:implementation}
discusses implementation, temporal accuracy, and related methods.
Numerical experiments are presented in Section~\ref{sec:numerics},
followed by concluding remarks in Section~\ref{sec:conclusion}.

\section{Preliminaries and a high-order linear predictor}\label{sec:prelim}

We introduce some notation and recall the energy balances. We then
present a generalized BDF predictor on a nonuniform time grid.

\subsection{Notation and continuous energy balances}

We denote by $\ip{\cdot}{\cdot}$ and $\norm{\cdot}$ the inner product
and norm in $L^2(\Omega)$, and use the same notation for scalar and vector
function spaces. Set
\[
 \bm V=\{\bv\in H_0^1(\Omega):\nabla\cdot\bv=0\}.
\]
There is a Poincar\'e constant $c_P>0$ such that
\begin{equation}\label{eq:poincare}
 \norm{\nabla\bv}^2\ge c_P\norm{\bv}^2,
 \qquad \bv\in H_0^1(\Omega).
\end{equation}
For periodic boundary conditions, we assume that the forcing and initial
velocity histories have zero mean, apply
\[
 \Pi_0\bv:=\bv-|\Omega|^{-1}\int_\Omega\bv\,dx
\]
to the right-hand side of each velocity predictor, and take the pressure
in $H^1_{\mathrm{per}}(\Omega)/\R$. The scalar corrections preserve zero
mean, so \eqref{eq:poincare} holds for the periodic iterates.

Taking the inner product of \eqref{eq:NS-mom} with $\bu$ gives
\begin{equation}\label{eq:continuous-L2}
 \frac{d}{dt}E(\bu)+\nu\norm{\nabla\bu}^2=\ip{\bfv}{\bu},
 \qquad E(\bv):=\frac12\norm{\bv}^2.
\end{equation}
In two dimensions, taking the inner product with $-\Delta\bu$ gives
\begin{equation}\label{eq:continuous-H1}
 \frac{d}{dt}K(\bu)=-D_2(\bu,p,\bfv),
 \qquad K(\bv):=\frac12\norm{\nabla\bv}^2,
\end{equation}
where
\begin{equation}\label{eq:D2}
 D_2(\bv,q,\bfv):=\nu\norm{\Delta\bv}^2+\ip{\bfv}{\Delta\bv}
 -\ip{\nabla q}{\Delta\bv}-\ip{(\bv\cdot\nabla)\bv}{\Delta\bv}.
\end{equation}

We introduce positive shifts $C_E,C_K>0$ and define
\begin{equation}\label{eq:shifted-energies}
 \cE(\bv):=E(\bv)+C_E,
 \qquad
 \cK(\bv):=K(\bv)+C_K.
\end{equation}
The shifts keep the energy ratios below well defined without changing
the continuous energy balances.

 Unless otherwise stated, the
stability constants may depend on $\Omega$, $\nu$, the correction degree
$k$, the forcing bound, and the fixed algorithmic parameters, but not on
the final time or the time grid. Dependence on the initial data is stated
separately.

\subsection{A flexible high-order predictor step}\label{sec:predictor}

Let
\begin{equation}\label{eq:time-grid}
 0=t_0<t_1<\cdots,\qquad \tau_{n+1}=t_{n+1}-t_n>0.
\end{equation}
For a fixed $T>0$, take $t_{N_T}=T$; for absorption, assume
$t_n\to\infty$. On a constant grid, $\tau_{n+1}=\dt$.

The scalar estimates in late sections apply to any well-defined linear predictor
with newly computed provisional velocity. We only need to make the following assumption on the linear predictor.

{\bf Assumption P.}  For any $n>0$, the linear predictor is of order $k$ for some $k\ge 1$ and produces provisional velocity  $\bar{\bm u}^{n+1} \in H^2(\Omega)$.

Examples include linearly implicit multistep and
Runge--Kutta methods.
As a particular example, we use the generalized BDF formulas of \cite{HuangShen2024} with extrapolated nonlinear term and pressure as the predictor. More precisely, for an order
$k\ge2$ and a shift $\beta\ge1$, define
\[
 \begin{aligned}
 A_{k,n}^{\beta}(\phi^{n+1})
  &=\sum_{q=0}^{k}a_{k,q}^{n}(\beta)\phi^{n+1-k+q},\\
 B_{k,n}^{\beta}(\phi^{n+1})
  &=\sum_{q=0}^{k-1}b_{k,q}^{n}(\beta)\phi^{n+2-k+q},\\
 C_{k,n}^{\beta}(\phi^n)
  &=\sum_{q=0}^{k-1}c_{k,q}^{n}(\beta)\phi^{n+1-k+q}.
 \end{aligned}
\]
Polynomial interpolation at the time nodes determines the coefficients.
The operator $A_{k,n}^{\beta}/\tau_{n+1}$ approximates the derivative
at $t_n+\beta\tau_{n+1}$, while $B_{k,n}^{\beta}$ and $C_{k,n}^{\beta}$
approximate the solution there to order $k$. On a constant grid these
reduce to the formulas in \cite{HuangShen2024}; $\beta=1$ gives classical BDF.
With $\bfv^{n+\beta}:=\bfv(t_n+\beta\tau_{n+1})$, set
\begin{equation}\label{eq:gbdf-predictor}
 \frac{A_{k,n}^{\beta}(\overline{\bu}^{n+1})}{\tau_{n+1}}
 -\nu\Delta B_{k,n}^{\beta}(\overline{\bu}^{n+1})
 +\big(C_{k,n}^{\beta}(\bu^n)\cdot\nabla\big)C_{k,n}^{\beta}(\bu^n)
 +\nabla C_{k,n}^{\beta}(p^n)=\bfv^{n+\beta}.
\end{equation}
Only $\overline{\bu}^{n+1}$ is unknown. To verify solvability, set
$t_*:=t_n+\beta\tau_{n+1}$. The leading Lagrange basis
functions for the derivative and diffusion interpolants are
\[
L_A(t)=\prod_{i=n+1-k}^{n}\frac{t-t_i}{t_{n+1}-t_i},
\qquad
L_B(t)=\prod_{i=n+2-k}^{n}\frac{t-t_i}{t_{n+1}-t_i}.
\]
Since $t_*\ge t_{n+1}>t_i$, all factors are positive and
\[
a_{k,k}^n(\beta)
=\tau_{n+1}L_A(t_*)\sum_{i=n+1-k}^{n}\frac1{t_*-t_i}>0,
\qquad b_{k,k-1}^n(\beta)=L_B(t_*)>0.
\]
Hence \eqref{eq:gbdf-predictor} has the form
$(aI-b\Delta)\overline{\bu}^{n+1}=\bg$ with $a,b>0$, and is uniquely
solvable. For an $L^2$ right-hand side, elliptic regularity gives
$\overline{\bu}^{n+1}\in H^2(\Omega)$ under the stated domain assumptions.

We use a well-defined linear starting procedure and apply the SAV
updates at every starting step. The derivative and diffusion terms use
the provisional-velocity histories, while the explicit convection uses
the corrected histories. High-order convergence requires order-compatible
starting values and the stability conditions of the chosen predictor.

\section{A one-SAV scheme with a uniform \texorpdfstring{$L^2$}{L2} bound}\label{sec:one-sav}

We construct a one-SAV scheme and prove a uniform $L^2$ bound and a
weighted $H^1$ dissipation estimate for the corrected velocity.

\subsection{Shifted \texorpdfstring{$L^2$}{L2} auxiliary variable}

Set $r(t)=\cE(\bu(t))$.  For parameters $\alpha_1,\beta_1>0$, define
\begin{equation}\label{eq:Q1-cont}
 Q_1(\bv,\bg):=\nu\norm{\nabla\bv}^2-\ip{\bg}{\bv}
 -\beta_1E(\bv)+\alpha_1\norm{\bg}^2.
\end{equation}
The scalar equation
\begin{equation}\label{eq:r-cont}
 r_t+\beta_1r=-\frac{r}{\cE(\bu)}Q_1(\bu,\bfv)
 +\alpha_1\norm{\bfv}^2+\beta_1C_E
\end{equation}
is equivalent to \eqref{eq:continuous-L2} when $r=\cE(\bu)$. We shall
choose $\alpha_1$ and $\beta_1$ so that $Q_1$ is nonnegative. We propose
the following scheme.

\textbf{Algorithm 1.} Given $\bu^n,p^n,r^n$ and the required history, perform the following steps.

\paragraph{Step 1: velocity predictor.}
Compute $\overline{\bu}^{n+1}$ from
\begin{equation}\label{eq:alg1-predictor}
 \frac{A_{k,n}^{\beta}(\overline{\bu}^{n+1})}{\tau_{n+1}}
 -\nu\Delta B_{k,n}^{\beta}(\overline{\bu}^{n+1})
 +\big(C_{k,n}^{\beta}(\bu^n)\cdot\nabla\big)C_{k,n}^{\beta}(\bu^n)
 +\nabla C_{k,n}^{\beta}(p^n)
 =\bfv^{n+\beta},
\end{equation}
with the prescribed boundary condition. Other linear predictors from
Section~\ref{sec:predictor} may replace \eqref{eq:alg1-predictor}.

\paragraph{Step 2: shifted SAV correction.}
Evaluate
\begin{equation}\label{eq:Q1-disc}
 Q_1^{n+1}:=\nu\norm{\nabla\overline{\bu}^{n+1}}^2
 -\ip{\bfv^{n+1}}{\overline{\bu}^{n+1}}
 -\beta_1E(\overline{\bu}^{n+1})
 +\alpha_1\norm{\bfv^{n+1}}^2,
\end{equation}
and compute $r^{n+1}$ from
\begin{equation}\label{eq:r-disc}
 \frac{r^{n+1}-r^n}{\tau_{n+1}}+\beta_1r^{n+1}
 =-\frac{r^{n+1}}{\cE(\overline{\bu}^{n+1})}Q_1^{n+1}
 +\alpha_1\norm{\bfv^{n+1}}^2+\beta_1C_E,
\end{equation}
or, equivalently,
\begin{equation}\label{eq:r-explicit}
 r^{n+1}=
 \frac{r^n+\tau_{n+1}\big(\alpha_1\norm{\bfv^{n+1}}^2+\beta_1C_E\big)}
 {1+\tau_{n+1}\beta_1+\tau_{n+1} Q_1^{n+1}/\cE(\overline{\bu}^{n+1})}.
\end{equation}
Define
\begin{equation}\label{eq:first-scaling}
 \xi^{n+1}:=\frac{r^{n+1}}{\cE(\overline{\bu}^{n+1})},
 \qquad
 \eta_k^{n+1}:=1-(1-\xi^{n+1})^k,
 \qquad
 \widetilde{\bu}^{n+1}:=\eta_k^{n+1}\overline{\bu}^{n+1}.
\end{equation}

\paragraph{Step 3: pressure recovery.}
Find $p^{n+1}\in H^1(\Omega)/\R$ such that
\begin{equation}\label{eq:pressure}
 \ip{\nabla p^{n+1}}{\nabla q}
 =\ip{\bfv^{n+1}-(\widetilde{\bu}^{n+1}\cdot\nabla)\widetilde{\bu}^{n+1}
 -\nu\curl\curl\widetilde{\bu}^{n+1}}{\nabla q},
 \qquad \forall q\in H^1(\Omega)/\R.
\end{equation}
Equation~\eqref{eq:pressure} is a constant-coefficient Poisson problem.
Set $\bu^{n+1}=\widetilde{\bu}^{n+1}$.

\begin{remark}[Order-preserving correction]\label{rem:eta-order}
On a fixed interval $[0,T]$, let
$h_{\max}=\max_{1\le j\le N_T}\tau_j$. The polynomial
$\eta_k(x)=1-(1-x)^k$ satisfies $\eta_k(1)=1$ and
$\eta_k(1+O(h_{\max}))=1+O(h_{\max}^k)$. Thus the correction is compatible
with an order-$k$ predictor, provided that $\xi^{n+1}=1+O(h_{\max})$.
\end{remark}

\subsection{Uniform \texorpdfstring{$L^2$}{L2} stability}

Assume
\begin{equation}\label{eq:param1}
 \alpha_1\ge\frac{2}{\nu c_P},
 \qquad
 0<\beta_1\le\frac{\nu c_P}{4}.
\end{equation}

We shall use the following damped scalar recursion
(cf.\ \cite[Lemma~1]{Shen1990}).

\begin{lemma}[Damped scalar recursion]\label{lem:shen-discrete}
Let $\tau_j>0$, $\lambda>0$, and let $a^j,b^j$ be nonnegative with
$b^j\le b$. If
\begin{equation}\label{eq:shen-recursion}
 \frac{a^{j+1}-a^j}{\tau_{j+1}}+\lambda a^{j+1}\le b^j,
\end{equation}
then, with $\Theta_n(\lambda):=\prod_{j=1}^{n}(1+\lambda\tau_j)^{-1}$
and $\Theta_0(\lambda)=1$,
\begin{equation}\label{eq:shen-bound}
 a^n\le\Theta_n(\lambda)a^0+
       \frac b\lambda\bigl(1-\Theta_n(\lambda)\bigr).
\end{equation}
Moreover, $\Theta_n(\lambda)\le(1+\lambda t_n)^{-1}$.
\end{lemma}
\begin{proof}
Set $q_j=(1+\lambda\tau_j)^{-1}$. Then
$a^{j+1}\le q_{j+1}a^j+(1-q_{j+1})b/\lambda$.
Subtract $b/\lambda$ and iterate to obtain \eqref{eq:shen-bound}.
The last assertion follows from
$\prod_{j=1}^{n}(1+\lambda\tau_j)\ge1+\lambda\sum_{j=1}^{n}\tau_j$.
\end{proof}

Define
\begin{equation}\label{eq:forcing-radius}
 B_1:=\alpha_1F_\infty^2+\beta_1C_E,\qquad
 R_\infty:=B_1/\beta_1,\qquad C_R:=\max\{r^0,R_\infty\}.
\end{equation}
For $R>0$, define the nondecreasing function
\begin{equation}\label{eq:Ceta}
 C_{\eta,k}(R):=\max_{0\le x\le R/C_E}x|P_{k-1}(x)|^2,\qquad
 P_{k-1}(x):=\sum_{\ell=0}^{k-1}(1-x)^\ell.
\end{equation}

The following estimates also hold for the first correction in
Algorithm~2.

\begin{theorem}[All-time $L^2$ bounds and velocity absorption]\label{thm:L2}
Let $d=2$ or $3$ and assume \eqref{eq:forcing}, \eqref{eq:param1}, and the predictor satisfies 
{\bf Assumption P}. Initialize
$r^0=\cE(\bu^0)>0$. Then, on any positive time grid, Algorithm~1 satisfies
\begin{equation}\label{eq:r-bound}
 0<r^n\le R_\infty+(r^0-R_\infty)_+\Theta_n(\beta_1)\le C_R,\qquad n\ge0,
\end{equation}
where $x_+:=\max\{x,0\}$. Furthermore,
\begin{equation}\label{eq:L2-final}
 \sup_{n\ge0}\norm{\bu^n}
 \le\max\{\norm{\bu^0},[2C_{\eta,k}(C_R)C_R]^{1/2}\}
 =:M_{L^2}.
\end{equation}
If $r^j\le R$ for $j\ge m$, then for every $n>m$,
\begin{equation}\label{eq:window-K}
 \begin{aligned}
 \sum_{j=m}^{n-1}\tau_{j+1}K(\widetilde{\bu}^{j+1})
 &\le \frac{C_{\eta,k}(R)}{\nu}\{r^m+B_1(t_n-t_m)\}
 \le A_R(t_n-t_m),\\
 A_R(L)&:=\frac{C_{\eta,k}(R)}{\nu}(R+B_1L).
 \end{aligned}
\end{equation}
If $t_n\to\infty$, then for each $\varepsilon>0$ the velocity eventually
lies in the $L^2$ ball with radius
\begin{equation}\label{eq:L2-absorbing-radius}
 M_{L^2,\mathrm{abs}}(\varepsilon)
 :=[2C_{\eta,k}(R_\infty+\varepsilon)
              (R_\infty+\varepsilon)]^{1/2}.
\end{equation}
The radius is independent of the initial data and time grid. Entry is
uniform over bounded sets of initial energies.
\end{theorem}

\begin{proof}
Young's inequality gives
\[
 |(\bfv^{n+1},\overline{\bu}^{n+1})|
 \le\frac{\nu c_P}{8}\norm{\overline{\bu}^{n+1}}^2
       +\frac{2}{\nu c_P}\norm{\bfv^{n+1}}^2.
\]
Using \eqref{eq:param1} and the Poincar\'e inequality yields
\begin{equation}\label{eq:Q1-lower}
 Q_1^{n+1}\ge
 \frac{\nu}{2}\norm{\nabla\overline{\bu}^{n+1}}^2
 +\frac{\nu c_P}{4}\norm{\overline{\bu}^{n+1}}^2\ge0.
\end{equation}
The denominator in \eqref{eq:r-explicit} is at least
$1+\beta_1\tau_{n+1}$, so $r^0>0$ implies $r^n>0$ by induction.
It follows from \eqref{eq:r-disc} that
\begin{equation}\label{eq:r-gronwall-recursion}
 \frac{r^{n+1}-r^n}{\tau_{n+1}}+\beta_1r^{n+1}
 \le\alpha_1\norm{\bfv^{n+1}}^2+\beta_1C_E\le B_1.
\end{equation}
Lemma~\ref{lem:shen-discrete} proves \eqref{eq:r-bound}.

To estimate the corrected velocity, we use the identity
\begin{equation}\label{eq:eta-polynomial}
 \eta_k(x)=1-(1-x)^k=xP_{k-1}(x)
\end{equation}
which implies $|\eta_k(x)|^2\le C_{\eta,k}(R)x$ for $0\le x\le R/C_E$.
Whenever $r^{n+1}\le R$,
\begin{equation}\label{eq:tilde-L2}
 \norm{\widetilde{\bu}^{n+1}}^2
 \le C_{\eta,k}(R)\frac{r^{n+1}}{\cE(\overline{\bu}^{n+1})}
                  \norm{\overline{\bu}^{n+1}}^2
 \le2C_{\eta,k}(R)r^{n+1}.
\end{equation}
Taking $R=C_R$ proves \eqref{eq:L2-final}. For the pressure estimates
below, set
\[
 M_0:=\max\{\norm{\bu^0},[2C_{\eta,k}(C_R)C_R]^{1/2}\}.
\]
Multiplying \eqref{eq:r-disc} at step $j$ by $\tau_{j+1}$, summing over
$j=m,\ldots,n-1$, and using \eqref{eq:Q1-lower}, we find
\[
 \frac{\nu}{2}\sum_{j=m}^{n-1}\tau_{j+1}\xi^{j+1}
          \norm{\nabla\overline{\bu}^{j+1}}^2
 \le r^m+B_1(t_n-t_m).
\]
Since
\begin{equation}\label{eq:K-eta-bound}
 K(\widetilde{\bu}^{j+1})
 \le\frac{C_{\eta,k}(R)}2\,\xi^{j+1}
                   \norm{\nabla\overline{\bu}^{j+1}}^2,
\end{equation}
we obtain \eqref{eq:window-K}.

For the absorbing bound, let $R_a=R_\infty+\varepsilon$.
Lemma~\ref{lem:shen-discrete} gives
\begin{equation}\label{eq:r-transient}
 r^n\le R_\infty+\frac{(r^0-R_\infty)_+}{1+\beta_1t_n}.
\end{equation}
Choose an index $n_a$ such that
\begin{equation}\label{eq:entry-index}
 t_{n_a}\ge\frac1{\beta_1}
       \left(\frac{(r^0-R_\infty)_+}{\varepsilon}-1\right)_+.
\end{equation}
Then $r^n\le R_a$ for every $n\ge n_a$, and \eqref{eq:tilde-L2}
with $R=R_a$ proves the eventual bound for $n\ge\max\{1,n_a\}$.
Replacing $r^0$ in \eqref{eq:entry-index} by a bound over an initial
set makes the entry criterion uniform over that set.
\end{proof}

\begin{remark}[Transient bounds and absorbing radii]\label{rem:absorption}
The all-time bound $M_{L^2}$ depends on the initial energy through $C_R$;
the absorbing radius \eqref{eq:L2-absorbing-radius} is independent of it
(cf.\ \cite{HumphriesStuart1994,JauberteauRosierTemam1990}).
If $\tau_j\le\tau_{\max}<\infty$, then $\log(1+x)\ge x/(1+x)$ gives
\[
 \Theta_n(\beta_1)
 \le\exp\!\left(-\frac{\beta_1t_n}{1+\beta_1\tau_{\max}}\right).
\]
Thus the initial contribution decays exponentially when the steps are
bounded, while \eqref{eq:r-transient} holds on arbitrary positive grids.
\end{remark}
\section{A two-SAV scheme with a uniform \texorpdfstring{$H^1$}{H1} bound in two dimensions}\label{sec:two-sav}

In this section, we assume that $d=2$ and the domain has a $C^3$
boundary. We introduce a second SAV and establish a uniform-in-time
$H^1$ bound by combining the weighted dissipation estimate for the
first SAV with the Stokes-pressure estimate of \cite{LiuLiuPego},
two-dimensional interpolation inequalities, and a discrete uniform
Gronwall argument. The scheme uses the velocity and pressure computed
in Algorithm~1 and requires no additional Helmholtz or Poisson solve.

\subsection{Shifted \texorpdfstring{$H^1$}{H1} auxiliary variable}

Let $s(t)=\cK(\bu(t))$ and choose $\alpha_2,\beta_2,\gamma_2>0$.  Define
\begin{equation}\label{eq:G2H2}
 G_2(\bv):=2\beta_2+2\gamma_2\norm{\nabla\bv}^2,
 \qquad
 H_2(\bfv):=\alpha_2\norm{\bfv}^2,
\end{equation}
so that
\[
 G_2(\bv)K(\bv)=\beta_2\norm{\nabla\bv}^2
 +\gamma_2\norm{\nabla\bv}^4.
\]
Using \eqref{eq:D2}, set
\begin{align}
 Q_2(\bv,q,\bfv)&:=D_2(\bv,q,\bfv)+G_2(\bv)K(\bv)+H_2(\bfv),\label{eq:Q2-cont}\\
 \widehat Q_2(\bv,q,\bfv)&:=Q_2(\bv,q,\bfv)+C_KG_2(\bv).\label{eq:Q2hat-cont}
\end{align}
Then
\begin{equation}\label{eq:s-cont}
 s_t=-\frac{s}{\cK(\bu)}\widehat Q_2(\bu,p,\bfv)+G_2(\bu)s+H_2(\bfv)
\end{equation}
reduces to \eqref{eq:continuous-H1} when $s=\cK(\bu)$, since the
additional terms cancel and the right-hand side equals
$-D_2(\bu,p,\bfv)$. Sufficient parameter conditions for a nonnegative
discrete residual will be given below.

\subsection{Algorithm 2: two SAVs}

Given $\bu^n,p^n,r^n,s^n$ and the required histories, perform the
following steps.

\paragraph{Steps 1--3.}
Compute the velocity predictor, first SAV correction, and pressure as
in Algorithm~1.

\paragraph{Step 4: recursive higher-energy update.}
Define
\begin{align}
 \cK_C^{n+1}&:=\cK(\widetilde{\bu}^{n+1})
 =\frac12\norm{\nabla\widetilde{\bu}^{n+1}}^2+C_K,\label{eq:KC-disc}\\
 G_2^{n+1}&:=2\beta_2+2\gamma_2\norm{\nabla\widetilde{\bu}^{n+1}}^2,\label{eq:G2-disc}\\
 H_2^{n+1}&:=\alpha_2\norm{\bfv^{n+1}}^2,\label{eq:H2-disc}
\end{align}
and
\begin{align}
 Q_2^{n+1}:={}&\nu\norm{\Delta\widetilde{\bu}^{n+1}}^2
 +\ip{\bfv^{n+1}}{\Delta\widetilde{\bu}^{n+1}}
 -\ip{\nabla p^{n+1}}{\Delta\widetilde{\bu}^{n+1}}\notag\\
 &-\ip{(\widetilde{\bu}^{n+1}\cdot\nabla)\widetilde{\bu}^{n+1}}
 {\Delta\widetilde{\bu}^{n+1}}
 +\beta_2\norm{\nabla\widetilde{\bu}^{n+1}}^2
 +\gamma_2\norm{\nabla\widetilde{\bu}^{n+1}}^4
 +\alpha_2\norm{\bfv^{n+1}}^2,\label{eq:Q2-disc}\\
 \widehat Q_2^{n+1}:={}&Q_2^{n+1}+C_KG_2^{n+1}.
 \label{eq:Q2hat-disc}
\end{align}

Compute a provisional scalar $\widehat S^{n+1}$ from
\begin{equation}\label{eq:Shat-eq}
 \frac{\widehat S^{n+1}-s^n}{\tau_{n+1}}
 =-\frac{\widehat S^{n+1}}{\cK_C^{n+1}}\widehat Q_2^{n+1}
 +G_2^{n+1}s^n+H_2^{n+1},
\end{equation}
or, explicitly,
\begin{equation}\label{eq:Shat}
 \widehat S^{n+1}=
 \frac{(1+\tau_{n+1} G_2^{n+1})s^n+\tau_{n+1} H_2^{n+1}}
 {1+\tau_{n+1}\widehat Q_2^{n+1}/\cK_C^{n+1}}.
\end{equation}

\paragraph{Step 5: scalar relaxation and final correction.}
Set
\begin{equation}\label{eq:relaxation}
 s^{n+1}:=\min\{\widehat S^{n+1},\cK_C^{n+1}\},
\end{equation}
and define
\begin{equation}\label{eq:second-scaling}
 \zeta^{n+1}:=\frac{s^{n+1}}{\cK_C^{n+1}},
 \qquad
 \lambda_k^{n+1}:=1-(1-\zeta^{n+1})^k,
 \qquad
 \bu^{n+1}:=\lambda_k^{n+1}\widetilde{\bu}^{n+1}.
\end{equation}

The relaxation \eqref{eq:relaxation} bounds the new scalar by the
shifted energy of $\widetilde{\bu}^{n+1}$ and thus transfers the
time-averaged estimate from the first SAV to $s^n$. Together with the
recursion \eqref{eq:Shat-eq}, this allows us to apply a uniform
Gronwall argument. The pressure is recovered before the second
rescaling so that the residual \eqref{eq:Q2-disc} uses the pair
$(\widetilde{\bu}^{n+1},p^{n+1})$.

\subsection{Positivity of the higher-energy residual}

\begin{lemma}[Two-dimensional pressure and convection estimates]\label{lem:pressure-2d}
Let $\bv\in H^2(\Omega)\cap H_0^1(\Omega)$ and let
$q\in H^1(\Omega)/\R$ satisfy
\begin{equation}\label{eq:pressure-abstract}
 \ip{\nabla q}{\nabla\varphi}
 =\ip{\bfv-(\bv\cdot\nabla)\bv-\nu\curl\curl\bv}{\nabla\varphi},
 \qquad \forall\varphi\in H^1(\Omega).
\end{equation}
Assume $\norm{\bv}\le M_0$. For every sufficiently small $\varepsilon>0$,
there exist constants $C_f,C_1,C_2,C_3,C_4>0$, independent of $\bv$,
$\bfv$, and the time grid, such that
\begin{align}
 \left|\ip{\nabla q}{\Delta\bv}\right|
 &\le\nu\left(\frac34+\frac32\varepsilon\right)\norm{\Delta\bv}^2
 +C_f\norm{\bfv}^2+C_1\norm{\nabla\bv}^2
 +C_2M_0^2\norm{\nabla\bv}^4,\label{eq:p-est}\\
 \left|\ip{(\bv\cdot\nabla)\bv}{\Delta\bv}\right|
 &\le\nu\epsilon\norm{\Delta\bv}^2
 +C_4M_0^2\norm{\nabla\bv}^4,\label{eq:nl-est}\\
 \left|\ip{\bfv}{\Delta\bv}\right|
 &\le\nu\epsilon\norm{\Delta\bv}^2+C_3\norm{\bfv}^2.
 \label{eq:f-est}
\end{align}
The constants may depend on $\Omega$, $\nu$, and $\varepsilon$.
\end{lemma}

\begin{proof}
Let $\mathbb P$ be the $L^2$ Leray projection and set
$\mathbb Q=I-\mathbb P$ and $N(\bv)=(\bv\cdot\nabla)\bv$.
It follows from \eqref{eq:pressure-abstract} that
\[
 \nabla q=\mathbb Q(\bfv-N(\bv))+\nu\nabla p_S,\qquad
 \nabla p_S=\mathbb Q(\Delta\bv-\nabla\nabla\cdot\bv).
\]
The Stokes-pressure estimate in \cite{LiuLiuPego} gives
\[
 \norm{\nabla p_S}^2
 \le(1/2+\varepsilon)\norm{\Delta\bv}^2+C_\varepsilon\norm{\nabla\bv}^2.
\]
This estimate holds for $\bv\in H^2\cap H_0^1$ without a
divergence-free constraint. Hence,
\[
 \nu|(\nabla p_S,\Delta\bv)|
 \le\nu(3/4+\varepsilon/2)\norm{\Delta\bv}^2
       +C_1\norm{\nabla\bv}^2.
\]
By the two-dimensional interpolation inequality and elliptic
regularity, we have
\[
 \norm{N(\bv)}
 \le C\norm{\bv}^{1/2}\norm{\nabla\bv}\norm{\Delta\bv}^{1/2}.
\]
Applying Young's inequality, we obtain, for every $\theta>0$,
\[
 |(N(\bv),\Delta\bv)|
 \le\nu\theta\norm{\Delta\bv}^2
    +C_{\nu,\theta}M_0^2\norm{\nabla\bv}^4.
\]
Since $\mathbb Q$ is an $L^2$ contraction, the same bound holds for
$\mathbb QN(\bv)$. Taking $\theta=\varepsilon/2$ and using
$|(\mathbb Q\bfv,\Delta\bv)|\le
(\nu\epsilon/2)\norm{\Delta\bv}^2+C_f\norm{\bfv}^2$, we obtain
\eqref{eq:p-est}. Taking $\theta=\varepsilon$ gives \eqref{eq:nl-est},
while \eqref{eq:f-est} follows from Young's inequality.
\end{proof}

Lemma~\ref{lem:pressure-2d} yields the following sufficient conditions
for positivity of the higher-energy residual.

\begin{proposition}[Positivity of the second residual]\label{prop:Q2-positive}
Let $0<\varepsilon<1/14$ and fix $\delta>0$. With the constants from
Lemma~\ref{lem:pressure-2d}, assume
\begin{equation}\label{eq:param2}
 \alpha_2\ge C_f+C_3+\delta,
 \qquad
 \beta_2\ge C_1+\delta,
 \qquad
 \gamma_2\ge(C_2+C_4)M_0^2+\delta.
\end{equation}
Then
\begin{equation}\label{eq:Q2-lower}
 Q_2^{n+1}\ge c_\Delta\norm{\Delta\widetilde{\bu}^{n+1}}^2
 +\delta\norm{\nabla\widetilde{\bu}^{n+1}}^2
 +\delta\norm{\nabla\widetilde{\bu}^{n+1}}^4
 +\delta\norm{\bfv^{n+1}}^2\ge0,
\end{equation}
where $c_\Delta=\nu(1/4-7\varepsilon/2)>0$. In particular,
$\widehat Q_2^{n+1}\ge0$.
\end{proposition}

\begin{proof}
The estimates \eqref{eq:p-est}--\eqref{eq:f-est} with
$\bv=\widetilde{\bu}^{n+1}$, $q=p^{n+1}$, and $\bfv=\bfv^{n+1}$
use the fraction
$3/4+3\varepsilon/2+\varepsilon+\varepsilon=3/4+7\varepsilon/2$ of the viscous
term. Together with \eqref{eq:param2}, this gives
\eqref{eq:Q2-lower}. Finally,
$\widehat Q_2^{n+1}=Q_2^{n+1}+C_KG_2^{n+1}\ge0$.
\end{proof}

\subsection{Weighted uniform Gronwall argument and the \texorpdfstring{$H^1$}{H1} bound}

We shall use the following variable-step version of the discrete
uniform Gronwall lemma \cite[Lemma~3]{Shen1990}. Its right-endpoint
weights agree with those in the dissipation estimate for the first SAV.

\begin{lemma}[Variable-step uniform Gronwall inequality]
\label{lem:uniform-gronwall}
Let $d^j,g^j,h^j$ be nonnegative and $\tau_j=t_j-t_{j-1}>0$. Suppose
\begin{equation}\label{eq:UG-recursion}
 d^{j+1}\le(1+\tau_{j+1}g^{j+1})d^j+\tau_{j+1}h^{j+1}.
\end{equation}
Fix integers $0\le m<n$, set $L=t_n-t_m$, and assume
\begin{equation}\label{eq:UG-windows}
 \sum_{j=m}^{n-1}\tau_{j+1}d^{j+1}\le A_d,\qquad
 \sum_{j=m}^{n-1}\tau_{j+1}g^{j+1}\le A_g,\qquad
 \sum_{j=m}^{n-1}\tau_{j+1}h^{j+1}\le A_h.
\end{equation}
Then
\begin{equation}\label{eq:UG-result}
 d^n\le\left(\frac{A_d}{L}+A_h\right)e^{A_g},
\end{equation}
without a restriction on adjacent step ratios. If $A_d,A_g,A_h$
are uniformly bounded over a family of intervals with $L\ge L_0>0$,
then \eqref{eq:UG-result} gives a uniform bound at their endpoints.
\end{lemma}

\begin{proof}
Iterating \eqref{eq:UG-recursion} from $i$ to $n$, we obtain
\[
 d^n\le d^i\prod_{\ell=i+1}^{n}(1+\tau_\ell g^\ell)
   +\sum_{\ell=i+1}^{n}\tau_\ell h^\ell
             \prod_{r=\ell+1}^{n}(1+\tau_r g^r).
\]
With empty products equal to one and empty sums equal to zero, this
also holds for $i=n$. Since the coefficients are nonnegative and
$1+x\le e^x$, we have
\[
 d^n\le e^{A_g}(d^i+A_h),\qquad m\le i\le n.
\]
Multiplying by $\tau_i$, summing over $i=m+1,\ldots,n$, and dividing
by $\sum_{i=m+1}^{n}\tau_i=L$, we obtain \eqref{eq:UG-result}
from \eqref{eq:UG-windows}. The last assertion follows directly.
\end{proof}

For Algorithm~2, the first SAV controls the right-endpoint average
$\sum\tau_{j+1}s^{j+1}$. We next establish the recursion required
by Lemma~\ref{lem:uniform-gronwall}.

\begin{lemma}[Recursion and upper energy bound]\label{lem:s-recursion}
Assume the conditions of Proposition~\ref{prop:Q2-positive} and
$s^0>0$. Then
\begin{equation}\label{eq:s-upper}
 0<s^{n+1}\le\cK_C^{n+1},
\end{equation}
and
\begin{equation}\label{eq:s-recursion}
 s^{n+1}\le(1+\tau_{n+1}G_2^{n+1})s^n+\tau_{n+1}H_2^{n+1}.
\end{equation}
\end{lemma}
\begin{proof}
By Proposition~\ref{prop:Q2-positive}, the denominator in
\eqref{eq:Shat} is at least one, and the numerator is positive whenever
$s^n>0$. Thus $\widehat S^{n+1}>0$, and \eqref{eq:relaxation}
gives \eqref{eq:s-upper} by induction. Dropping the denominator
in \eqref{eq:Shat} and using $s^{n+1}\le\widehat S^{n+1}$
gives \eqref{eq:s-recursion}.
\end{proof}

We now prove the uniform $H^1$ bound.

\begin{theorem}[All-time $H^1$ stability on arbitrary positive grids]
\label{thm:H1}
Let $T>0$, $0=t_0<\cdots<t_{N_T}=T$, and $d=2$. Assume
\eqref{eq:forcing}, \eqref{eq:param1}, \eqref{eq:param2},  and the predictor satisfies 
{\bf Assumption P}.
For some fixed \(\sigma\in(0,1]\), suppose that
\[
  \bu^0=\overline{\bu}^0\in H^{1+\sigma}(\Omega),
\]
and that the starting histories satisfy
\[
  \bu^j,\overline{\bu}^j\in H^2(\Omega),
  \qquad j=1,\ldots,k-1,
  \qquad
  p^j\in H^1(\Omega)/\mathbb R,
  \qquad j=0,\ldots,k-1.
\]
Require all velocity histories to satisfy the prescribed homogeneous
Dirichlet conditions. Initialize
\[
  r^0=\cE(\bu^0),
  \qquad
  s^0=\cK(\bu^0),
\]
and apply both SAV updates and the scalar relaxation at every starting
step. Prescribe the initial pressure \(p^0\) in
\(H^1(\Omega)/\mathbb R\) and apply pressure recovery to each newly
computed positive-time velocity. Let the subsequent values be generated
by Algorithm~2. Then
\begin{equation}\label{eq:H1-final}
 \norm{\bu^n}\le M_{L^2},\qquad
 \norm{\nabla\bu^n}\le M_{H^1},\qquad 0\le n\le N_T.
\end{equation}
The constants may depend on the initial data and fixed algorithmic
parameters, but are independent of $T$, $N_T$, and all step sizes.
The bounds require neither an upper step-size bound nor a restriction
on adjacent step ratios, and hold for all $n$ on an infinite
positive grid.
\end{theorem}

\begin{proof}
By Theorem~\ref{thm:L2}, the first correction satisfies
$\norm{\widetilde{\bu}^n}\le M_0$ and \eqref{eq:window-K}.
For $R>0$, set
\begin{equation}\label{eq:window-constants}
 D_R:=A_R(2)+2C_K,\qquad
 \Gamma_R:=4\beta_2+4\gamma_2A_R(2),\qquad
 H_*:=2\alpha_2F_\infty^2.
\end{equation}
Fix the reference window length at $1$, independently of the grid.
For any interval $[t_m,t_n]$ of length at most $2$ with
$r^j\le R$ for $j\ge m$, \eqref{eq:window-K} and
\eqref{eq:s-upper} imply
\begin{equation}\label{eq:weighted-averages}
 \sum_{j=m}^{n-1}\tau_{j+1}s^{j+1}\le D_R,\qquad
 \sum_{j=m}^{n-1}\tau_{j+1}G_2^{j+1}\le\Gamma_R,\qquad
 \sum_{j=m}^{n-1}\tau_{j+1}H_2^{j+1}\le H_*.
\end{equation}
Here $G_2^{j+1}=2\beta_2+4\gamma_2K(\widetilde{\bu}^{j+1})$.

Taking $R=C_R$, we shall prove
\begin{equation}\label{eq:s-global}
 s^n\le M_S:=[\max\{s^0,D_{C_R}\}+H_*]e^{\Gamma_{C_R}}
 \quad\text{for every }n.
\end{equation}
First, for any step $\tau_\ell\ge1$, the one-step estimate
\eqref{eq:window-K} gives
\[
 K(\widetilde{\bu}^{\ell})
 \le\frac{C_{\eta,k}(C_R)}{\nu}
       \left(\frac{C_R}{\tau_\ell}+B_1\right)
 \le A_{C_R}(1).
\]
Thus $s^\ell\le A_{C_R}(1)+C_K\le D_{C_R}$.

For a fixed endpoint $n$, let $\ell$ be the index of the last step
of size at least $1$ among steps $1,\ldots,n$, and set $\ell=0$
if there is no such step. All steps after $\ell$ are less than $1$.
If $t_n-t_\ell<1$, iterating \eqref{eq:s-recursion} and using
\eqref{eq:weighted-averages}, we obtain
\[
 s^n\le e^{\Gamma_{C_R}}(s^\ell+H_*)
 \le e^{\Gamma_{C_R}}(\max\{s^0,D_{C_R}\}+H_*).
\]
The case $\ell=n$ is included.
If $t_n-t_\ell\ge1$, choose the largest $m\ge\ell$ such that
$t_n-t_m\ge1$. Since $\tau_{m+1}<1$, we have
$1\le t_n-t_m<2$. Lemma~\ref{lem:uniform-gronwall}, with
$d^j=s^j$, $g^j=G_2^j$, $h^j=H_2^j$, and
\eqref{eq:weighted-averages}, yields
\[
 s^n\le(D_{C_R}+H_*)e^{\Gamma_{C_R}}\le M_S.
\]
This proves \eqref{eq:s-global}.

To estimate the velocity, note that $0<\zeta^{n+1}\le1$ and
\[
 0\le\lambda_k(z)=1-(1-z)^k\le\min\{1,kz\},\qquad 0\le z\le1.
\]
It follows that $\lambda_k(z)^2\le\lambda_k(z)\le kz$, and hence
\begin{equation}\label{eq:H1-transfer}
 K(\bu^{n+1})
 \le k\zeta^{n+1}K(\widetilde{\bu}^{n+1})
 \le ks^{n+1}\le kM_S.
\end{equation}
Also $\norm{\bu^{n+1}}\le\norm{\widetilde{\bu}^{n+1}}\le M_0$.
Thus \eqref{eq:H1-final} holds with
$M_{L^2}=M_0$ and $M_{H^1}=(2kM_S)^{1/2}$, including $n=0$
because $M_S\ge s^0\ge K(\bu^0)$.
\end{proof}

\begin{remark}[Regularity of the starting values]
At the first full multistep update, $\bu^0$ enters the derivative
and explicit convection, but not the diffusion interpolation. In two
dimensions, we have
\[
  \|(w\cdot\nabla)v\|_{L^2}
  \le \|w\|_{L^\infty}\|\nabla v\|_{L^2}
  \le C_\sigma
       \|w\|_{H^{1+\sigma}}\|v\|_{H^{1+\sigma}},
  \qquad \sigma>0.
\]
Thus its convection contribution belongs to $L^2$ for
$H^{1+\sigma}$ initial data. The estimate fails when both factors
are only in $H^1$, so the argument does not cover $\sigma=0$.

For the generalized diffusion interpolation, each independent old
provisional velocity with a nonzero diffusion coefficient must satisfy
$\Delta\overline{\bu}^j\in L^2$, or equivalently
$\overline{\bu}^j\in H^2$ under the stated boundary and domain
assumptions. Together with $p^j\in H^1/\R$, this gives an $L^2$
right-hand side for the Helmholtz equation. Elliptic regularity yields
$\overline{\bu}^k\in H^2$, scalar rescaling preserves this regularity,
and pressure recovery gives $p^k\in H^1/\R$. Repeating the argument
gives the same regularity at all subsequent steps.

These sufficient regularity conditions can be weakened if a diffusion
coefficient vanishes or a history combination has additional regularity.
For classical BDF, diffusion is evaluated only at the new time level,
so $H^{1+\sigma}$ starting velocity histories suffice for the first
full multistep predictor. The $H^2$ regularity justifies each solve
and residual evaluation; no uniform-in-time $H^2$ bound is assumed
or obtained.
\end{remark}

Applying the preceding argument after the first scalar enters its
absorbing interval gives the following result.

\begin{corollary}[Eventual $H^1$ velocity bounds]\label{cor:H1-absorption}
Fix the algorithmic parameters and consider initialized histories
satisfying the hypotheses of Theorem~\ref{thm:H1} with these
parameters. Let $t_n\to\infty$. For $\varepsilon>0$, set
$R_a=R_\infty+\varepsilon$ and choose $n_a$ by \eqref{eq:entry-index}.
Then
\begin{equation}\label{eq:s-absorbing}
 s^n\le (D_{R_a}+H_*)e^{\Gamma_{R_a}}
 =:M_{S,\mathrm{abs}}(\varepsilon)
 \qquad\text{whenever }t_n-t_{n_a}\ge1.
\end{equation}
Consequently, at all such endpoints,
\begin{equation}\label{eq:H1-absorbing}
 \norm{\bu^n}_{H^1}^2
 \le 2C_{\eta,k}(R_a)R_a+2kM_{S,\mathrm{abs}}(\varepsilon).
\end{equation}
Thus \eqref{eq:H1-absorbing} defines an absorbing ball for the
admissible initial class, with radius independent of the particular
history and the time grid. The entry criterion is uniform over
bounded sets of admissible initial energies.
\end{corollary}

\begin{proof}
By Theorem~\ref{thm:L2}, $r^j\le R_a$ for $j\ge n_a$.
Hence \eqref{eq:window-K} and \eqref{eq:weighted-averages} hold with
$R=R_a$ on every interval starting at or after $n_a$.
We apply the proof of Theorem~\ref{thm:H1} starting at $n_a$.
If a step of size at least $1$ occurs among steps
$n_a+1,\ldots,n$, let $\ell$ be the index of the last such step.
Then $s^\ell\le D_{R_a}$.
For $t_n-t_\ell<1$, direct iteration gives the bound
$(D_{R_a}+H_*)e^{\Gamma_{R_a}}$. Otherwise, the remaining steps
are less than $1$, and we apply Lemma~\ref{lem:uniform-gronwall}
on an interval $[t_m,t_n]$ with $1\le t_n-t_m<2$.
If no such large step occurs after $n_a$, the condition
$t_n-t_{n_a}\ge1$ gives the same backward interval.
This proves \eqref{eq:s-absorbing} without a bound on $s^{n_a}$.
Finally, \eqref{eq:tilde-L2} and \eqref{eq:H1-transfer} yield
\eqref{eq:H1-absorbing}.
\end{proof}

\begin{remark}[Parameter dependence and the meaning of absorption]
\label{rem:parameter-absorption}
To impose \eqref{eq:param2} uniformly on an initial class whose
all-time $L^2$ bound satisfies $M_0\le M_*$, fix
$\gamma_2\ge(C_2+C_4)M_*^2+\delta$.
Corollary~\ref{cor:H1-absorption} then gives a common absorbing radius
for this class. The radius depends on the fixed $\gamma_2$, but not
on the particular initial history or its $H^1$ norm. Enlarging the
class may require a different $\gamma_2$ and hence a different radius.
The all-time bound retains its dependence on $r^0$ and $s^0$ through
$C_R$ and $M_S$.
\end{remark}

\section{Implementation and comparison}\label{sec:implementation}

In this section, we discuss the implementation and temporal accuracy of
the two schemes, and compare their long-time stability properties and
computational costs with those of existing methods.

\subsection{Implementation and computational cost}

After spatial discretization, the generalized BDF predictor takes the form
\begin{equation}\label{eq:matrix-helmholtz}
 \left(\frac{a_{k,k}^{n}(\beta)}{\tau_{n+1}}M
       +\nu b_{k,k-1}^{n}(\beta)A\right)\overline U^{n+1}=R^n,
\end{equation}
where $M$ and $A$ are the fixed mass and Laplace matrices, respectively,
and $R^n$ contains only known quantities. The pressure is recovered by
solving the Poisson equation
\begin{equation}\label{eq:matrix-poisson}
 A_pP^{n+1}=G(\widetilde U^{n+1},F^{n+1}).
\end{equation}
For a fixed spatial mesh and a constant time step, we can factorize the
matrices in \eqref{eq:matrix-helmholtz}--\eqref{eq:matrix-poisson} once
and reuse the factorizations at subsequent steps. For variable time
steps, the spatial operators remain unchanged, but the scalar weights
in the Helmholtz matrix may require a new factorization. Finite-element,
finite-difference, and spectral discretizations on nonperiodic domains
preserve this constant-coefficient structure in space.

In addition to the field solves and the evaluation of residuals and
norms, Algorithm~1 requires the scalar division \eqref{eq:r-explicit}
and one velocity rescaling. Algorithm~2 adds the division
\eqref{eq:Shat}, the minimum in \eqref{eq:relaxation}, and a second
rescaling. Hence, the stronger $H^1$ bound requires no additional
Helmholtz or Poisson solve. The coefficients of the linear systems are
independent of the numerical velocity, and no nonlinear iteration is
needed. These costs are for the multistep schemes presented above; the
cost of a Runge--Kutta predictor depends on its stages.

The stability results are established for the time-discrete schemes.
To obtain mesh-uniform estimates for a fully discrete scheme, one needs
discrete counterparts of the Poincar\'e, pressure, and interpolation
estimates used in the analysis. In particular, for standard $C^0$ finite
elements, replacing the continuous Laplacian in the higher-energy
residual by an elementwise Laplacian requires further analysis.

\subsection{Temporal accuracy}

The correction polynomials convert a first-order error in the SAV ratios
into higher-order perturbations of the velocity. On a fixed time
interval $[0,T]$, let $h_{\max}=\max_{1\le j\le N_T}\tau_j$. Then
\begin{equation}\label{eq:order-map}
 \xi^n=1+O(h_{\max})\ \Longrightarrow\
 \eta_k^n=1+O(h_{\max}^k),\qquad
 \zeta^n=1+O(h_{\max})\ \Longrightarrow\
 \lambda_k^n=1+O(h_{\max}^k).
\end{equation}
For the minimum relaxation, we have
\[
 1-\zeta^{n+1}
 =\left(1-\frac{\widehat S^{n+1}}{\cK_C^{n+1}}\right)_+,
\]
so the minimum does not increase the provisional ratio error relative
to the current shifted energy. The ratio estimates in
\eqref{eq:order-map} require an error analysis of the coupled schemes;
they do not follow from uniform boundedness alone.

For the BDF scheme, we retain the provisional histories in the
derivative and diffusion terms and use corrected velocities in the
explicit terms, as in the degree-$k$ correction of
\cite{HuangShen2021}. For a one-step method restarted from the corrected
velocity, a perturbation of order $h^k$ at each step can accumulate to
order $h^{k-1}$ on a constant grid. In this case, a degree-$(k+1)$
correction is one possible choice, provided that the corresponding
ratio-consistency estimates hold. The scalar stability analysis applies
to both choices and, more generally, to any fixed correction degree.

High-order convergence further requires sufficiently smooth solutions
and forcing, starting values of the appropriate order, and the stability
and time-grid conditions of the predictor. Thus, boundedness on arbitrary
positive grids does not imply accuracy on such grids or zero-stability
of classical BDF formulas at all orders. A finite-time error analysis of
the coupled schemes remains to be developed.

\begin{remark}[Relaxation variants]
Any alternative relaxation satisfying
\[
 0<s^{n+1}\le\cK_C^{n+1},\qquad s^{n+1}\le\widehat S^{n+1}
\]
preserves the two inequalities used in the stability proof. Accuracy
also requires consistency of the scalar ratio. For example, one half of
the minimum satisfies both inequalities but need not converge to the
continuous energy. Hence, these stability conditions alone do not ensure
the accuracy of an alternative relaxation.
\end{remark}

\subsection{Comparison with existing long-time stable methods}

In Table~\ref{tab:comparison}, we compare methods for which
uniform-in-time bounds have been proved. The third column gives the
order and time-step conditions of the cited long-time analysis. For the
present schemes, we list the correction degree, since the stability
analysis for an arbitrary degree does not establish the corresponding
convergence order.

\begin{table}[htbp]
\caption{Comparison of numerical schemes with proved uniform-in-time
bounds for the Navier--Stokes equations. Here ``unconditional'' means
that the cited long-time estimate requires no CFL-type or
data-dependent time-step restriction.}
\label{tab:comparison}
\centering
\scriptsize
\setlength{\tabcolsep}{3pt}
\renewcommand{\arraystretch}{1.12}
\begin{tabularx}{\textwidth}{@{}>{\raggedright\arraybackslash}p{3.0cm}
 >{\raggedright\arraybackslash}p{4.0cm}
 >{\raggedright\arraybackslash}p{2.8cm}>{\raggedright\arraybackslash}X@{}}
\toprule
Method and reference & Setting and uniform bound & Order and long-time condition
& Principal solves per step \\
\midrule
Fully implicit Euler and Crank--Nicolson
\cite{ToneEuler,ToneCN}
& No-slip, 2D; velocity $L^2$ and $H^1$
& Orders $1$ and $2$; the $H^1$ analysis requires a time-step condition
  (CFL type for Crank--Nicolson)
& Nonlinear velocity--pressure problem \\

Linear backward Euler \cite{Ju}
& Primitive variables, 2D; no-slip or periodic; bounds in $V$ and $D(A)$,
  where $A$ is the Stokes operator
& Order $1$; bounds uniform with respect to $\dt$
& Linear Oseen/Stokes problem with a convection coefficient depending on the solution \\

Linearly extrapolated BDF2
\cite{AkbasKayaRebholz,RebholzTone,Shiue2025}
& Primitive variables, no-slip; velocity $L^2$ in 2D/3D and $H^1$ in 2D
& Order $2$; $L^2$ unconditional; $H^1$ requires a mild time-step
  restriction in the cited analysis
& Linear Oseen/Stokes problem with a convection coefficient depending on the solution \\

Velocity--vorticity IMEX finite elements \cite{Heister}
& Periodic, 2D; velocity and vorticity in both $L^2$ and $H^1$
& BDF2 and a first-order variant; unconditional
& Two sequential linear Oseen-type field solves; the vorticity solve is decoupled \\

Explicit-convection periodic vorticity schemes
\cite{GottliebToneWang,Wang2012,ChengWang}
& Periodic, 2D; vorticity $L^2$ and $H^1$, with higher Sobolev bounds
  in the later schemes
& Orders $1$--$4$; time-step restriction
& One or two constant-coefficient Poisson/heat solves, depending on the scheme \\

Earlier high-order SAV--IMEX scheme \cite{HuangShen2021}
& Periodic, 2D/3D; velocity $H^1$ in 2D and $L^2$ in 3D;
  zero external force in the analyzed model
& BDF orders $1$--$5$; unconditional
& One constant-coefficient Poisson-type velocity solve and an explicit scalar update \\

Forced SAV--BDF2 scheme \cite{HanWang2024,LinShiue2025}
& Periodic, no-slip, or free-slip; velocity $L^2$ in 2D/3D;
  vorticity $L^2$ and $H^1$ in the 2D free-slip setting
& Order $2$; unconditional; the later analysis allows nonsmooth initial data
  and gives a $\dt$-independent energy bound
& Two constant-coefficient Brinkman/Stokes solves and one linear scalar update \\

Variable-step ETD--mr-SAV schemes
\cite{WangETD2026,WangLinearETD2026}
& Periodic, 2D; vorticity $L^2$
& Variable-step order $2$ (with an embedded first-order companion);
  unconditional
& \cite{WangETD2026}: two heat solves and a scalar cubic;
  \cite{WangLinearETD2026}: two heat solves and a linear scalar recovery by
  Laplace--Talbot inversion \\

Present one-SAV scheme
& Velocity $L^2$ for no-slip problems in 2D/3D
& Any correction degree $k$; arbitrary positive variable steps
& One vector Helmholtz solve, one Poisson solve, and one scalar division \\

Present two-SAV scheme
& No-slip, 2D; velocity $H^1$
& Any correction degree $k$; arbitrary positive variable steps
& The same field solves, plus one scalar division and a scalar relaxation \\
\bottomrule
\end{tabularx}
\end{table}

The methods in Table~\ref{tab:comparison} differ mainly in the systems
to be solved at each step. Fully implicit methods require nonlinear
solves. Linearly extrapolated BDF2 and velocity--vorticity methods lead
to linear Oseen-type systems whose coefficients depend on the numerical
solution. The periodic vorticity methods with explicit convection use
constant-coefficient operators, but their cited long-time estimates
require time-step restrictions.

Among the SAV methods, the high-order scheme in
\cite{HuangShen2021} treats the unforced periodic problem. The forced
SAV--BDF2 schemes in \cite{HanWang2024,LinShiue2025} provide
uniform-in-time bounds and require constant-coefficient Brinkman/Stokes
solves. The variable-step ETD methods in
\cite{WangETD2026,WangLinearETD2026} treat the two-dimensional periodic
vorticity equation and require either a scalar cubic solve or a linear
auxiliary recovery by numerical Laplace inversion.

Our schemes use the velocity--pressure formulation and require only
decoupled Helmholtz/Poisson solves and explicit scalar updates. The
second SAV strengthens the velocity bound in two dimensions without
additional field solves. Moreover, the scalar estimates are valid for
the admissible predictors and arbitrary positive time grids.

\FloatBarrier
\section{Numerical examples}\label{sec:numerics}

In this section, we present three groups of numerical tests in two
dimensions with no-slip boundary conditions. We first verify the
temporal convergence of the two-SAV scheme. We then compare the BDF3
predictor with and without SAV corrections, and finally examine the
efficiency of adaptive time stepping for a problem with a short initial
transient and persistent long-time motion. We use constant time steps
in the first two tests and variable time steps in the third. The
algorithmic parameters are fixed within each refinement study.

\subsection{Convergence test}
\label{subsec:num-convergence}

We take $\Omega=(-1,1)^2$, $\nu=1$, and $T=1$, and choose the forcing
such that the exact velocity and pressure are given by
\begin{equation}\label{eq:num-noslip-exact}
 \begin{aligned}
 \bm u(x,y,t)&=\pi
 \begin{pmatrix}
 \sin(2\pi y)\sin^2(\pi x)\sin(t-\tfrac12)\\
 -\sin(2\pi x)\sin^2(\pi y)\sin(t-\tfrac12)
 \end{pmatrix},\\
 p(x,y,t)&=\sin(t-\tfrac12)\cos(\pi x)\sin(\pi y).
 \end{aligned}
\end{equation}

We use the Legendre--Galerkin method in space with polynomial degree
$N_x=N_y=32$. In each coordinate, the velocity basis consists of
$L_j-L_{j+2}$, $j=0,\ldots,N_x-2$, where $L_j$ is the Legendre polynomial
of degree $j$. For the pressure, we use all polynomials up to degree $32$
and fix the constant mode by the zero-mean condition. We test the
generalized BDF schemes of orders $k=2,3,4$ with the shifts and time steps
\begin{equation}\label{eq:num-gbdf-parameters}
 (k,\beta)=(2,2),\ (3,3),\ (4,8),
 \qquad
 \dt=\frac1{80},\frac1{160},\frac1{320},\frac1{640}.
\end{equation}
For Algorithm~2, we set
\begin{equation}\label{eq:num-conv-sav-parameters}
 \begin{gathered}
 C_E=C_K=100,\qquad
 \alpha_1=\frac4{\pi^2},\qquad\beta_1=\frac{\pi^2}{8},\\
 \alpha_2=5,\qquad\beta_2=\gamma_2=10^{-2}.
 \end{gathered}
\end{equation}
We use correction degree $k$ and the minimum relaxation
\eqref{eq:relaxation}. The derivative and diffusion terms use the
provisional velocity histories, while the explicit convection uses the
final corrected velocities. The pressure is recovered from the first
corrected velocity before applying the second SAV correction.

\begin{figure}[!htbp]
	\centering
	\includegraphics[width=0.8\textwidth]{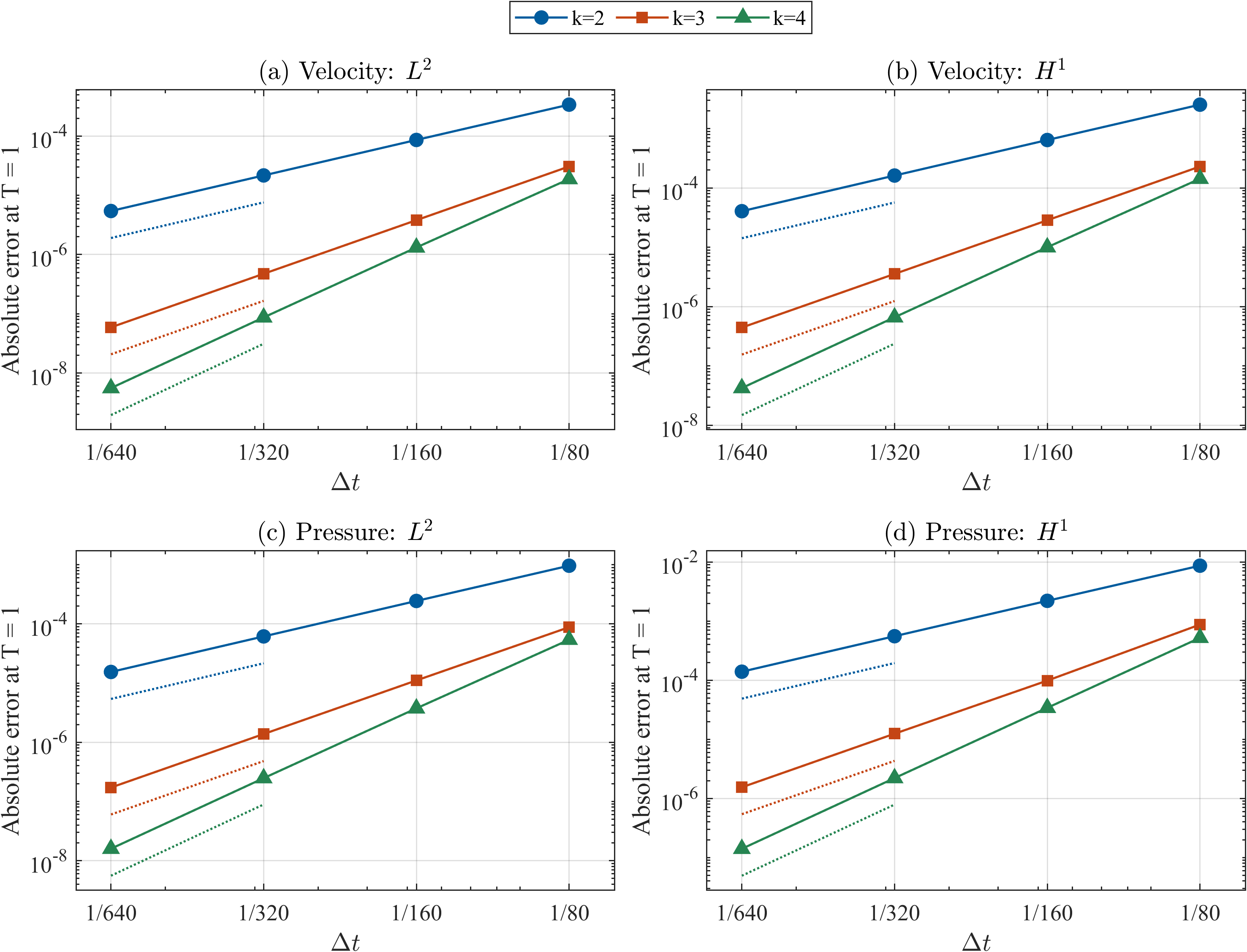}
	\caption{Temporal convergence of Algorithm~2 for
		\eqref{eq:num-noslip-exact} at $T=1$. Solid lines with filled markers
		show the results for $k=2,3,4$, with different colors and marker
		shapes. Dotted segments of the same colors indicate slopes
		$2$, $3$, and $4$, respectively.}
	\label{fig:num-convergence-fields}
\end{figure}

For $q=\bu,p$, we define the errors at the final time by
\begin{equation}\label{eq:num-errors}
 e_{q,X}(\dt)=\|q^{N_T}-q(T)\|_X,
 \quad X=L^2,H^1.
\end{equation}
Here $\|v\|_{H^1}^2=\|v\|^2+\|\nabla v\|^2$, and the numerical and
exact pressures have zero mean. From
Figure~\ref{fig:num-convergence-fields}, we observe the expected orders
$2$, $3$, and $4$ for both velocity and pressure in the two norms over
the tested range of time steps.

\begin{figure}[!htbp]
 \centering
 \includegraphics[width=0.8\textwidth]{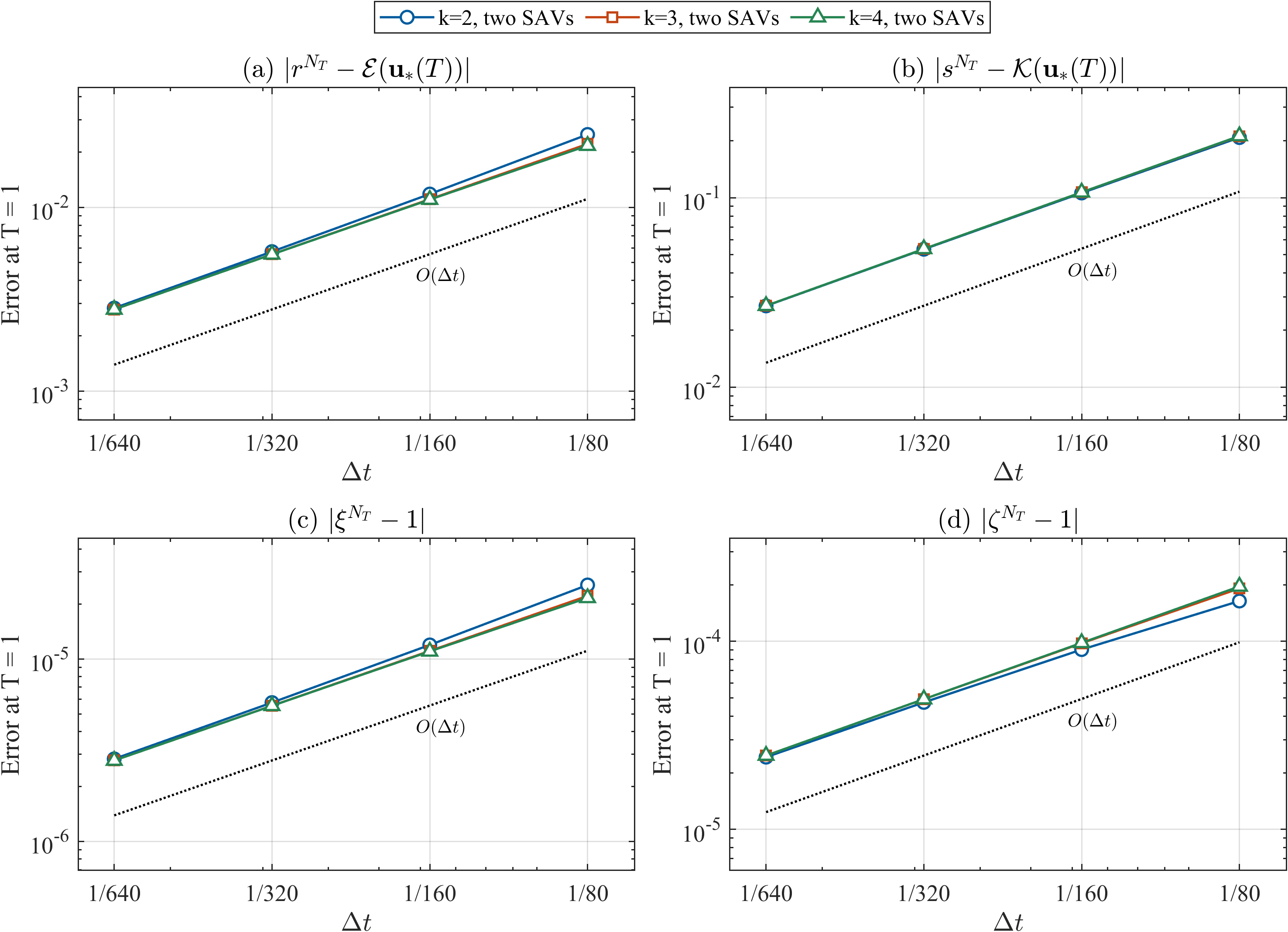}
 \caption{Scalar and ratio errors for Algorithm~2 in the convergence
 test. Dotted lines indicate slope one. The scalar errors are first
 order, while the velocity and pressure errors have the higher orders
 shown in Figure~\ref{fig:num-convergence-fields}.}
 \label{fig:num-convergence-scalars}
\end{figure}

Next, let $\bu_*$ denote the exact velocity in
\eqref{eq:num-noslip-exact}, and define the scalar and ratio errors by
\begin{equation}\label{eq:num-scalar-errors}
 \begin{aligned}
 e_r&=|r^{N_T}-\cE(\bu_*(T))|,&
 e_s&=|s^{N_T}-\cK(\bu_*(T))|,\\
 e_\xi&=|\xi^{N_T}-1|,&
 e_\zeta&=|\zeta^{N_T}-1|.
 \end{aligned}
\end{equation}
In Figure~\ref{fig:num-convergence-scalars}, we observe approximately
first-order convergence of these errors. Thus, in this test, first-order
scalar errors coexist with higher-order velocity and pressure
convergence, consistent with \eqref{eq:order-map}. The
minimum relaxation is active in these computations. All computed SAV
residuals are nonnegative, and all update denominators are positive.

\subsection{Effect of the SAV corrections on stability}
\label{subsec:num-stability}

We compare the BDF3 predictor with and without the two SAV corrections
on the unit disk $D=\{(x,y):x^2+y^2<1\}$ and the square $S=(-1,1)^2$,
both with no-slip boundary conditions. In each comparison, the predictor
and spatial discretization are unchanged, so the difference is due to
the scalar updates and rescalings.

We first consider the unit disk. Set $\nu=10^{-2}$, $A=1.05$,
$\sigma=x^2+y^2$, and prescribe
\begin{equation}\label{eq:num-df-exact}
\begin{aligned}
 \Psi&=(1-\sigma)^2\big[A+10^{-4}e^{-t}
                  \operatorname{Re}[(x+\mathrm i y)^7]\big],
 &\bu_*&=(\partial_y\Psi,-\partial_x\Psi),\\
 p_*&=A^2\big[8(\sigma-\sigma^2+\sigma^3/3)-2\big].
\end{aligned}
\end{equation}
The velocity is smooth and divergence free. We evaluate the bounded
force analytically from $\bfv=\partial_t\bu_*-\nu\Delta\bu_*+
(\bu_*\cdot\nabla)\bu_*+\nabla p_*$, independently of the numerical
solution.

We use the Fourier--Jacobi velocity space
$(1-\sigma)\mathbb P_{N-2}(D)^2$, with $N=14$ and $19\times65$ quadrature
points. The BDF3 formula uses the provisional velocity histories, while
convection uses $3\bu^n-3\bu^{n-1}+\bu^{n-2}$. We use the scalar updates,
pressure recovery, and minimum relaxation of Algorithm~2 with cubic
correction polynomials. We take $\dt=0.04$, $T=120$, and
\[
 C_E=C_K=1,\quad \alpha_1=4/(\nu\pi^2),\quad
 \beta_1=\nu\pi^2/8,\quad \alpha_2=50,\quad
 \beta_2=\gamma_2=10^{-2}.
\]
Both computations use the same forcing, spatial operators, and RK4
startup with $20$ substeps on each of the first two intervals. The
computation without SAV omits only the scalar updates and rescalings.
Neither computation uses filtering or resets to the exact solution.

\begin{figure}[!htbp]
\centering
\includegraphics[width=0.8\linewidth]{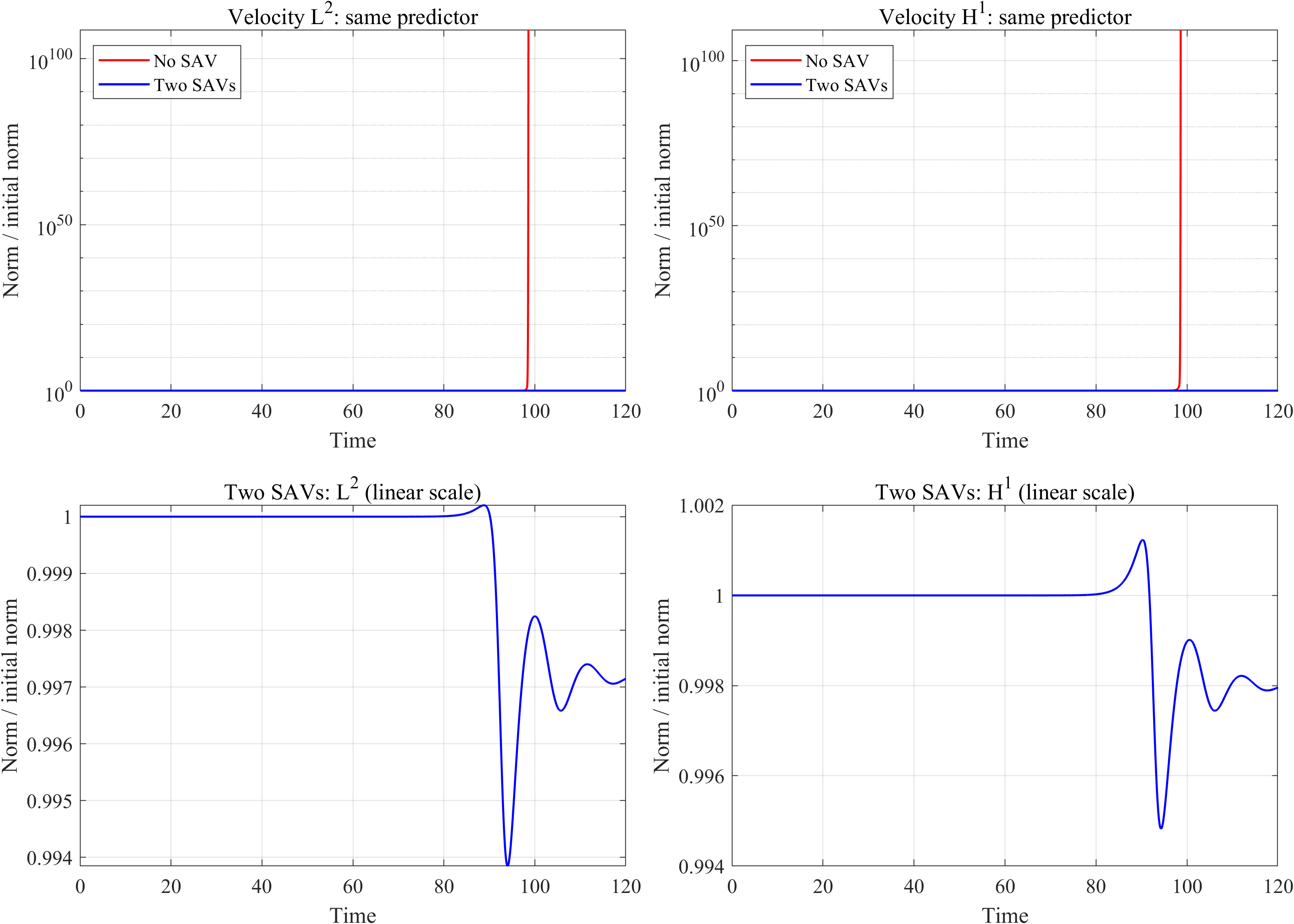}
\caption{Velocity norms normalized by their initial values ($N=14$,
$\dt=0.04$). Top: results with and without SAV on logarithmic scales.
Bottom: two-SAV results on linear scales.}
\label{fig:num-df-norms}
\end{figure}

\begin{figure}[!htbp]
\centering
\includegraphics[width=0.8\linewidth]{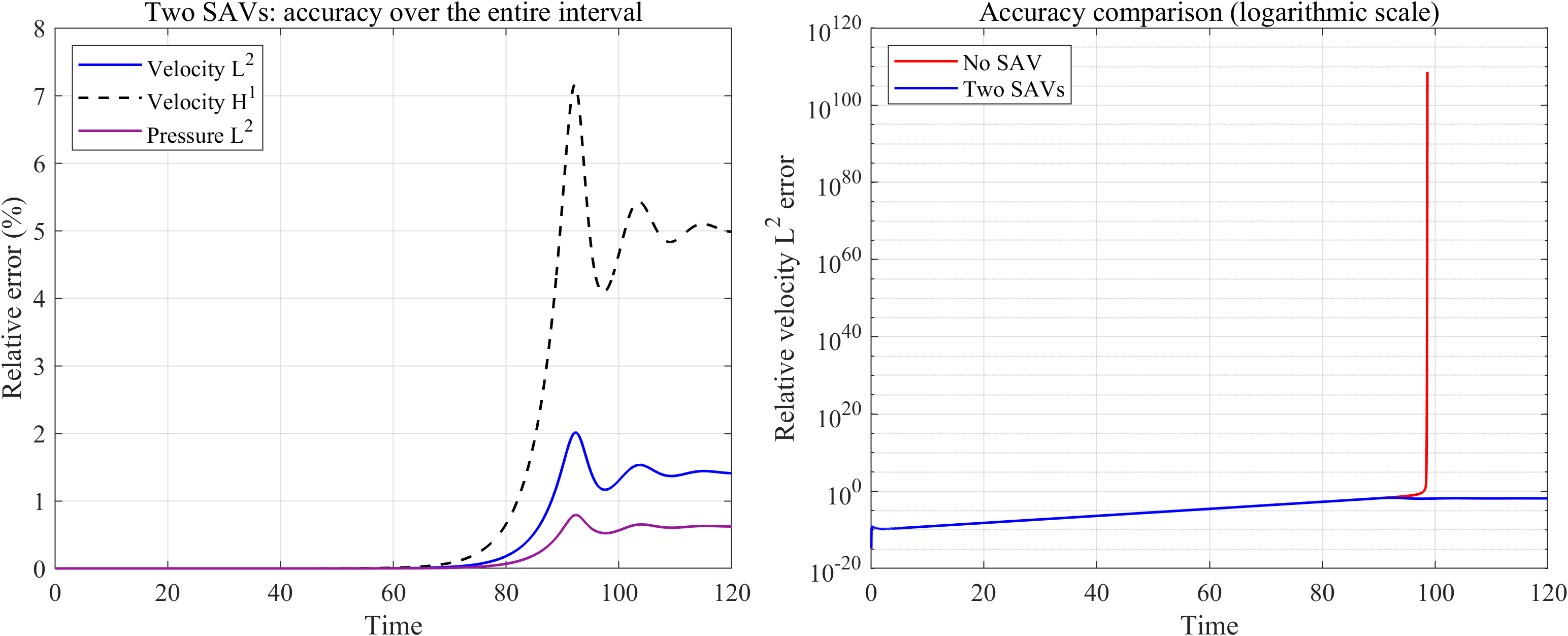}
\caption{Relative errors. Left: two-SAV velocity and pressure errors in
percent. Right: velocity $L^2$ errors with and without SAV on a
logarithmic scale.}
\label{fig:num-df-errors}
\end{figure}

In Figure~\ref{fig:num-df-norms}, we observe that the uncorrected
predictor becomes unstable near $t=99$, while the two-SAV solution
remains bounded through all $3000$ steps. In
Figure~\ref{fig:num-df-errors}, we compare the numerical solutions with
the exact solution. For the two-SAV scheme, the maximum relative errors
over the full interval are
$2.0118\%$ and $7.1616\%$ for velocity in $L^2$ and the full $H^1$ norm,
respectively, and $0.7921\%$ for pressure in $L^2$.

This example also shows that boundedness does not ensure accuracy. For
$N=18$ and $\dt=0.04$, the two-SAV solution remains bounded, but its
maximum relative velocity $L^2$ error is about $21\%$. With
$\dt=0.02$, both methods give accurate results for $N=14,18$. Thus,
although the scalar corrections control the velocity norms, the time
step must still be chosen to resolve the solution.

The \(H^1\) analysis in Section~\ref{sec:two-sav} uses a Stokes-pressure estimate for
domains with \(C^3\) boundaries. We next consider a square to test the
scheme on a polygonal domain, which is outside this assumption.

We consider $S=(-1,1)^2$ with homogeneous no-slip boundary
conditions. Let $D=\{x^2+y^2<1\}$ and retain the disk data
\[
 \psi_D=(1-\sigma)^2[A+\varepsilon e^{-t}\operatorname{Re}[(x+iy)^7]],
 \quad \bu_D=(\partial_y\psi_D,-\partial_x\psi_D),
 \quad \sigma=x^2+y^2,
\]
with $A=1.05$, $\varepsilon=10^{-4}$, and
$p_D=A^2[8(\sigma-\sigma^2+\sigma^3/3)-2]$.
Writing $\bfv_D=\partial_t\bu_D-\nu\Delta\bu_D+
(\bu_D\cdot\nabla)\bu_D+\nabla p_D$, we prescribe
\begin{equation}\label{eq:square-zero-extension}
 \bu_0=\mathbf{1}_D\bu_D(\cdot,0),\qquad
 \bfv(\cdot,t)=\mathbf{1}_D\bfv_D(\cdot,t)
 \quad\text{on }S.
\end{equation}
Since $\bu_D=0$ on $\partial D$, the extended initial velocity belongs
to $H_0^1(S)^2$ and is weakly divergence free. We do not assume that
the extended disk field is an exact time-dependent solution on the
square.

We use a Legendre--Galerkin BDF3 predictor in
$[(1-x^2)(1-y^2)\mathbb Q_{12,12}]^2$,
with $\nu=10^{-2}$, $\dt=0.04$, and $T=120$. Both computations use the
same projected initial data and RK4 starting histories, with 20 substeps
per initial interval. For the two-SAV scheme, we use cubic corrections
and the minimum relaxation, with
\[
 C_E=C_K=1,\quad \alpha_1=\frac4{\nu\pi^2},\quad
 \beta_1=\frac{\nu\pi^2}{8},\quad \alpha_2=50,\quad
 \beta_2=\gamma_2=10^{-2}.
\]
The computation without SAV omits only the scalar updates and
rescalings. We evaluate the initial-data and forcing integrals on their
disk support, and the velocity integrals on the entire square.

\begin{figure}[!htbp]
 \centering
 \includegraphics[width=\textwidth]{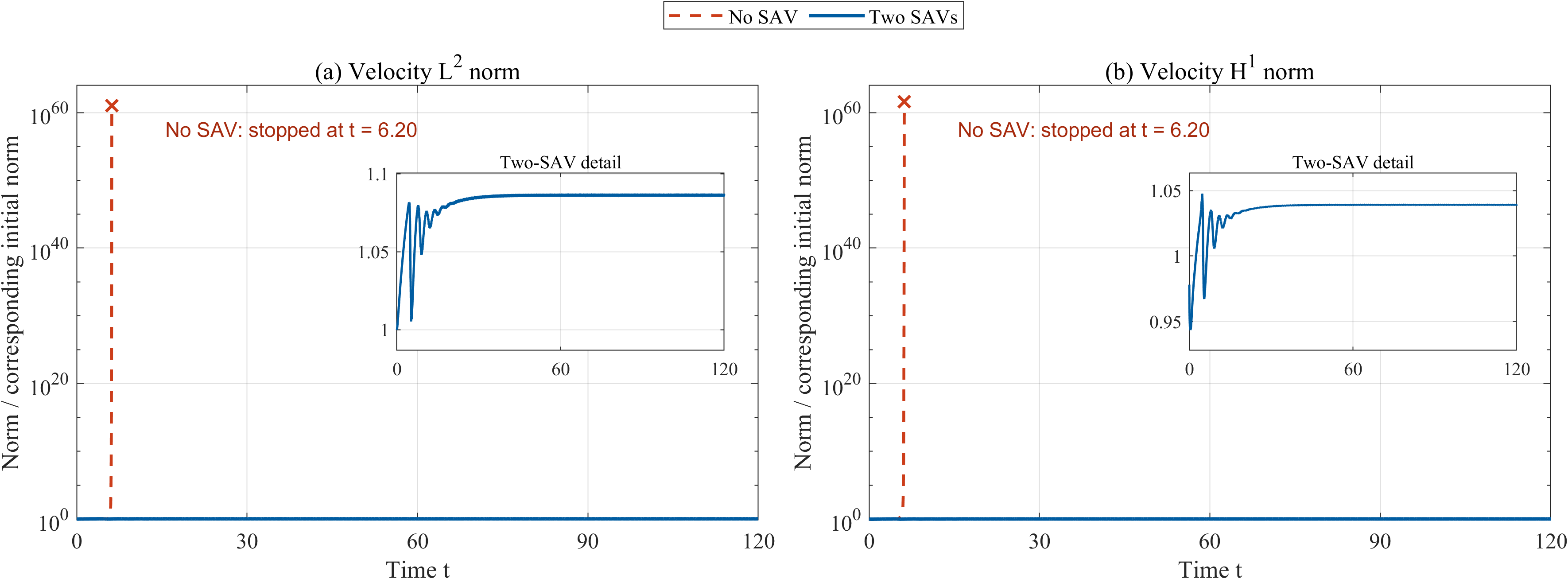}
 \caption{Stability comparison on the square at $\dt=0.04$. Velocity
 norms are normalized by the corresponding continuous initial norms.
 The main axes show all recorded values; the computation without SAV
 terminates at the runaway threshold before $T=120$. Insets show enlarged
 views of the two-SAV norm histories.}
 \label{fig:square-stability}
\end{figure}

In Figure~\ref{fig:square-stability}, the computation without SAV
becomes unstable:
$\|\bu^n\|/\|\bu_0\|$ exceeds $10^{60}$ at $t=6.20$.
The two-SAV computation reaches $T=120$, with maximum normalized $L^2$
and full $H^1$ norms of $1.0873$ and $1.0476$, respectively. All computed
SAV residuals and update denominators are positive. We observe the same
qualitative stability behavior as on the disk. Extending the stability
analysis to polygonal domains remains an open problem.

\subsection{Adaptive long-time integration}
\label{subsec:adaptive-long-time}

We now examine the efficiency of adaptive time stepping. We combine
the two-SAV corrections with the two-stage IMEX--SDIRK2 predictor and
embedded step-selection strategy of \cite{Liao2026}, replacing its SAV
equation by our scalar equations and minimum relaxation. Each step
starts from the corrected velocity. For this second-order predictor,
we use the cubic rescalings $\eta_3=1-(1-\xi)^3$ and
$\lambda_3=1-(1-\zeta)^3$, as discussed in
Section~\ref{sec:implementation}. We set the embedded relative velocity
tolerance to $5\times10^{-3}$ and the tolerance for the scalar indicator
$\max\{|1-\xi|,|1-\zeta|\}$ to $0.5$. The safety factor is $0.9$, the
initial step is $10^{-3}$, and the step bounds are $[10^{-7},0.2]$.
If a step is rejected, all trial velocity, pressure, and scalar values
are discarded.

On the no-slip unit disk, let $\sigma=x^2+y^2$, $A=1.05$, and prescribe
\begin{equation}\label{eq:adaptive-exact}
\begin{aligned}
 a(t)&=1+e^{-100t}+0.2\sin(0.1t),\\
 \Psi(t)&=(1-\sigma)^2\big[Aa(t)+10^{-4}e^{-t}
                    \operatorname{Re}[(x+\mathrm i y)^7]\big],
 &\bu_*&=(\partial_y\Psi,-\partial_x\Psi),\\
 p_*(t)&=a(t)^2A^2\big[8(\sigma-\sigma^2+\sigma^3/3)-2\big].
\end{aligned}
\end{equation}
We evaluate the bounded force analytically from
$\bfv=\partial_t\bu_*-\nu\Delta\bu_*+(\bu_*\cdot\nabla)\bu_*+\nabla p_*$.
We use $\nu=10^{-2}$, $N=14$, $19\times65$ quadrature points, and the
SAV parameters of the preceding disk example. The solution has an
initial layer on the time scale $0.01$, followed by persistent motion
with period $20\pi$. All computations use the same spatial operators,
forcing, and initial data, with final time $T=1000$ and no solution
resets.

\begin{figure}[!htbp]
\centering
\includegraphics[width=0.8\linewidth]{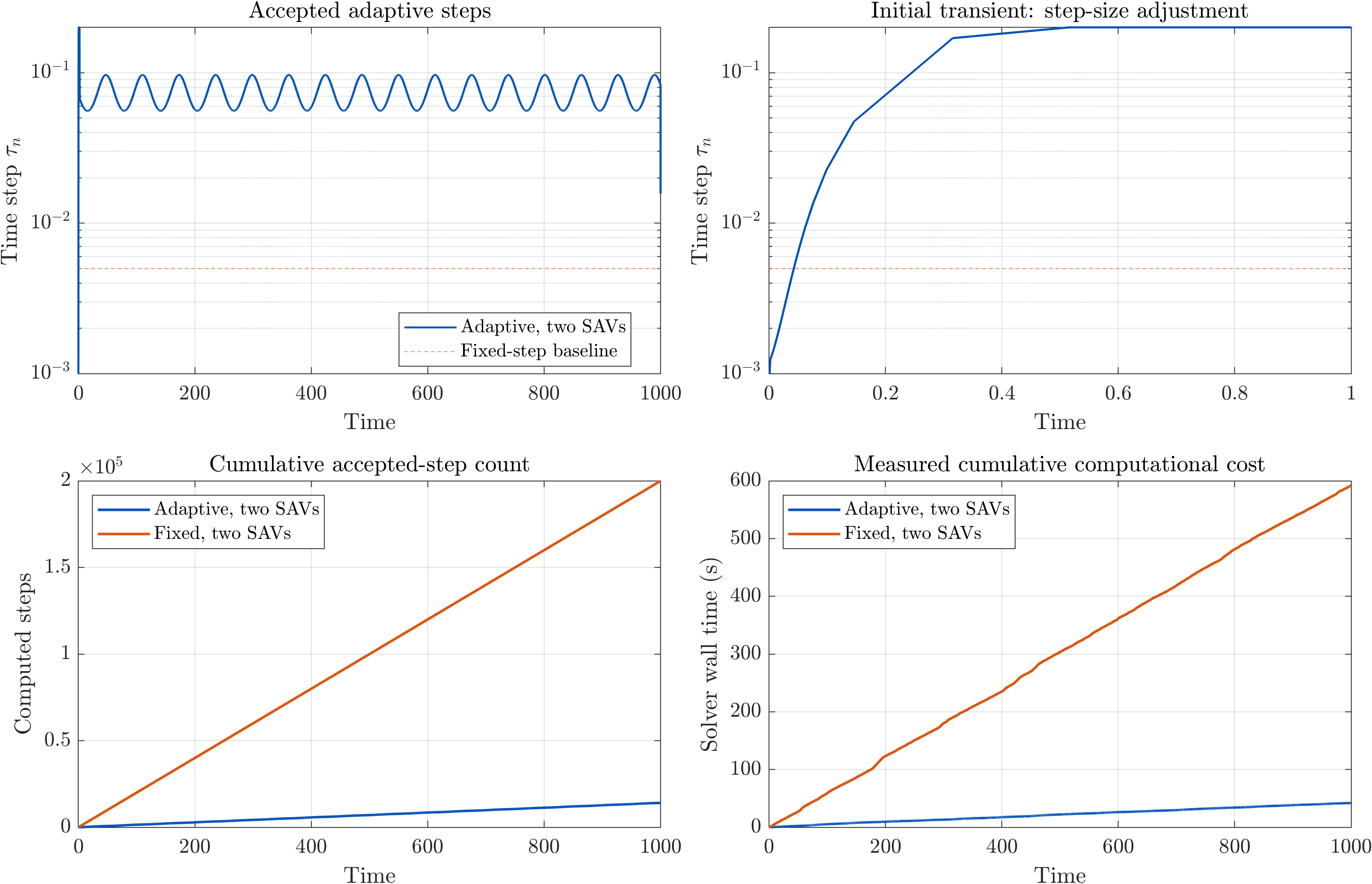}
\caption{Efficiency of the adaptive two-SAV scheme. Top: accepted step
sizes over the full interval and during the initial transient. Bottom:
cumulative step counts and solver wall times compared with the
fixed-step two-SAV scheme.}
\label{fig:adaptive-efficiency}
\end{figure}

\begin{figure}[!htbp]
\centering
\includegraphics[width=0.8\linewidth]{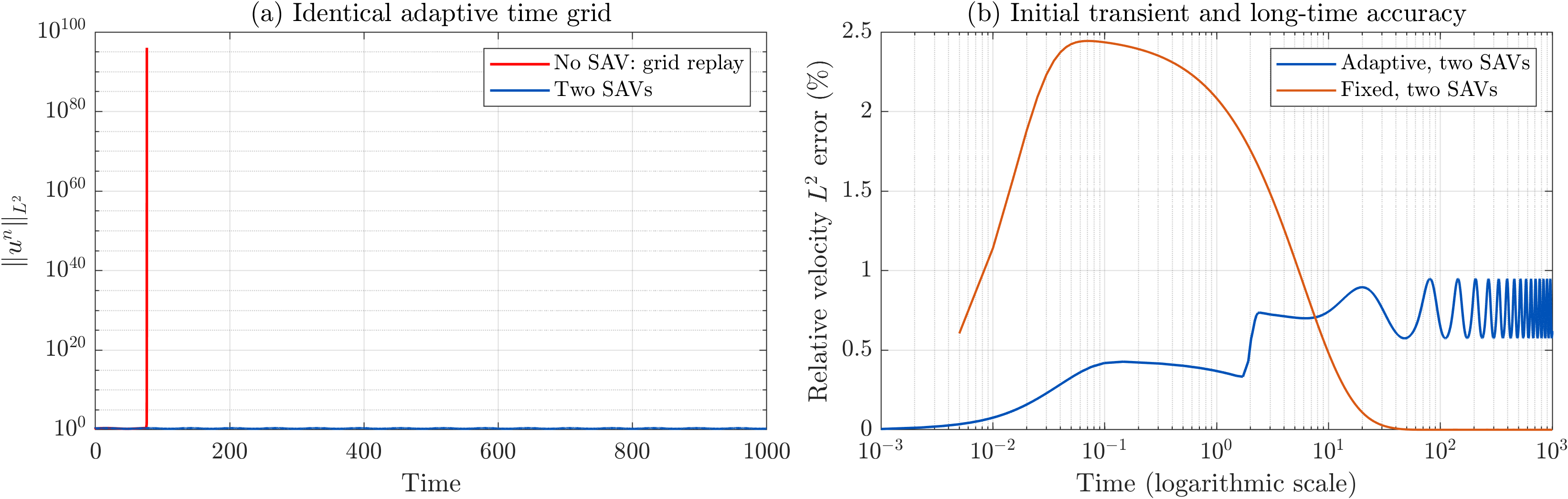}
\caption{(a) Velocity $L^2$ norms with and without SAV on the same
adaptive time grid; the computation without SAV terminates near $t=77$.
(b) Relative velocity $L^2$ errors of the adaptive and fixed-step two-SAV
schemes, with a logarithmic time axis.}
\label{fig:adaptive-stability}
\end{figure}

We first test several constant time steps up to $T=200$. Among these,
$\tau=0.005$ is the largest for which the maximum relative velocity
error is at most $5\%$ in $L^2$. Both the fixed-step and adaptive two-SAV
schemes satisfy this requirement over $[0,1000]$, with maximum errors
$2.444\%$ and $0.949\%$, respectively. In
Figure~\ref{fig:adaptive-efficiency}, we observe that the adaptive scheme
uses small steps to resolve the initial layer and larger steps
thereafter. The number of accepted steps decreases from $200000$ to
$14083$, with two rejected trials, and the reported solver wall time
decreases from $591.92$ to $42.01$ seconds, giving a speedup of $14.09$.
Both methods meet the same accuracy requirement, but their errors are
not matched. In particular, the fixed-step solution is more accurate
after the initial transient (Figure~\ref{fig:adaptive-stability}(b)).

Finally, we repeat the computation on the accepted adaptive grid without
SAV corrections. The uncorrected solution becomes unstable near $t=77$,
while the two-SAV solution remains bounded to $T=1000$
(Figure~\ref{fig:adaptive-stability}(a)), with a maximum relative
velocity $H^1$ error of $2.640\%$. A computation without SAV that selects
its own adaptive steps also reaches $T=1000$, with maximum relative
velocity $L^2$ error $0.737\%$. Thus, the comparison on the same grid
shows the stabilizing effect of the SAV corrections, while a suitable
adaptive grid also allows the uncorrected method to remain stable in
this example. The reported speedup compares adaptive and fixed steps
within the two-SAV method.

\section{Concluding remarks}\label{sec:conclusion}
We constructed in this paper two SAV correction steps for the forced
incompressible Navier--Stokes equations. A main feature is that their
formulation and long-time stability analysis are independent of the
particular linear predictor. The same correction steps apply to any
linear predictor satisfying Assumption~P. The first correction uses
a damped, shifted kinetic-energy variable to obtain uniform-in-time
velocity $L^2$ bounds in two and three dimensions. In two dimensions,
the second correction combines a higher-energy scalar equation with
a relaxation to yield a uniform-in-time $H^1$ bound. Under the stated
assumptions, the estimates hold on
arbitrary positive time grids. We also established eventual absorbing
bounds; for the two-SAV scheme, these apply to the admissible initial
class associated with fixed algorithmic parameters.

The analysis combines the weighted dissipation estimate from the
first SAV with the Stokes-pressure estimate and a discrete uniform
Gronwall argument. By retaining right-endpoint weights and treating
large steps directly, we avoid restrictions on step sizes and
adjacent step ratios. The scalar estimates do not require an
energy identity for the predictor. Moreover, the second SAV requires
no additional Helmholtz or Poisson solve.

The numerical results exhibit second-, third-, and fourth-order
convergence for the tested BDF schemes and illustrate the effect of
the scalar corrections on long-time stability. In the adaptive
example, variable steps reduced the reported computational time by
about a factor of fourteen relative to the fixed-step two-SAV
computation, while both satisfied the same accuracy requirement.
This result demonstrates the benefit of adaptivity for the tested
implementation.

The uniform-in-time bounds established here do not by themselves
imply uniform-in-time accuracy. Solutions of the forced
Navier--Stokes equations need not be stable with respect to
perturbations of the initial data. Consequently, small discretization
errors may be amplified over long times even when both the exact
and numerical solutions remain bounded. Error estimates of the
expected convergence order with constants independent of the final
time generally require additional stability assumptions on the
exact solution, such as uniform exponential stability
\cite{heywood1986finite}. A finite-time error analysis of the coupled
schemes and the approximation of long-time statistical quantities
remain to be studied.

\bibliographystyle{plain}
\bibliography{long_time_stability_NS}

\end{document}